\documentstyle[epsf,12pt]{article}
\def\Bbb R{{\rm \bf R}}
\def\proclaim#1{\vskip2mm{\bf #1}\em}
\def\endproclaim{\em \vskip2mm}
\def\tag#1{\eqno(#1)}
\def\gathered{\begin{array}{c}}
\def\endgathered{\end{array}}
\def\text{\mbox}

\begin{document}

\title {Extracting the geometry of an obstacle and a zeroth-order coefficient of a boundary condition
via the enclosure method using a single reflected wave
over a finite time interval}
\author{Masaru IKEHATA\footnote{
Laboratory of Mathematics,
Institute of Engineering,
Hiroshima University,
Higashi-Hiroshima 739-8527, JAPAN}}
%\date{}
\maketitle

\begin{abstract}
This paper considers an inverse problem for the classical wave equation in
an exterior domain.  It is a mathematical interpretation of an
inverse obstacle problem which employs the dynamical scattering
data of acoustic wave over a finite time interval.  It is assumed
that the wave satisfies a Robin type boundary condition with an
unknown variable coefficient.

The wave is generated by the
initial data localized outside the obstacle and observed over a
finite time interval at the same place as the support of the
initial data.  It is already known that, using the enclosure
method, one can extract the {\it maximum sphere} whose exterior
encloses the obstacle, from the data.

In this paper, it is shown
that the enclosure method enables us to extract also:

(i)  a quantity
which indicates the {\it deviation} of the geometry between the maximum
sphere and the boundary of the obstacle at the first reflection
points of the wave;

(ii)  the value of the coefficient of the boundary
condition at an arbitrary first reflection point of the wave provided,
for example, the surface of the obstacle is known in a
neighbourhood of the point.

Another new obtained knowledge is that: the enclosure method can cover
the case when the data are taken over a {\it sphere} whose centre coincides with that of the support
of an initial data and yields corresponding results to (i) and (ii).

\noindent
AMS: 35R30, 35L05, 35J05

\end{abstract}

%\tableofcontents

\section{Introduction}

Produce a {\it single wave} with {\it compact support}
at the initial state outside an unknown obstacle and observe the reflected wave at some place not far way 
from the obstacle over a {\it finite time interval}.  The poblem
of extracting information about the geometry and property of the surface of the obstacle
from this observed wave is a {\it proto-type} of so-called {\it inverse obstacle problems}
and the solution may have many applications to, e.g., {\it sonar}, {\it radar} imaging.
In this paper, we consider an inverse problem for the classical wave equation in
an exterior domain which is a mathematical interpretation of this inverse obstacle problem.

Let us describe a mathematical formulation of the problem.
Let $D$ be a nonempty bounded open subset of $\Bbb R^3$ with $C^2$-boundary
such that $\Bbb R^3\setminus\overline D$ is connected.
Let $0<T<\infty$.
Let $f\in L^2(\Bbb R^3)$ satisfy $\text{supp}\,f\cap\overline D=\emptyset$.
Let $u=u_f(x,t)$ denote the weak solution of the following initial boundary value problem for the classical wave equation:
$$\begin{array}{c}
\displaystyle
\partial_t^2u-\triangle u=0\,\,\text{in}\,(\Bbb R^3\setminus\overline D)\times\,]0,\,T[,\\
\\
\displaystyle
u(x,0)=0\,\,\text{in}\,\Bbb R^3\setminus\overline D,\\
\\
\displaystyle
\partial_tu(x,0)=f(x)\,\,\text{in}\,\Bbb R^3\setminus\overline D,\\
\\
\displaystyle
\frac{\partial u}{\partial\nu}-\gamma(x)\partial_tu-\beta(x)u=0\,\,\text{on}\,\partial D\times\,]0,\,T[.
\end{array}
\tag {1.1}
$$
Here $\nu$ denotes the unit {\it outward normal} to $D$ on $\partial D$.
The coefficients $\gamma(\ge 0)$ and $\beta$ in the boundary condition in (1.1) are essentially bounded functions on $\partial D$.
The weak solution
for the wave equation over a finite time interval is the same as used in \cite{IE0, IE3}
which comes from \cite{DL}.

%Stop exaggerated formulization.
%Set up a problem by the simplest possible formulization.
Given $p\in\Bbb R^3$ define $d_{\partial D}(p)=\inf_{x\in\partial D}\vert x-p\vert$.
We denote by $\Lambda_{\partial D}(p)$ the set of all points $q\in\partial D$
such that $\vert x-p\vert=d_{\partial D}(p)$.  We call this the {\it first reflector} from $p$ to $\partial D$
and the points in the first reflector are called the {\it first reflection points} from $p$ to $\partial D$.

In this paper, first we consider the following inverse problem.

{\bf Problem I.}
Let $B$ be a {\it known} open ball centred at $p\in\Bbb R^3$ and with radius $\eta$ such that $\overline B\cap\overline D=\emptyset$.
Let $\chi_B$ denote the characteristic function of $B$ and set $f=\chi_B$.

(i)  Assume that $D$ is {\it unknown} and both $\gamma$ and $\beta$ {\it unknown}.  Extract information
about the location and shape of $D$ from the data $u_f(x,t)$ given at all
$x\in B$ and $t\in\,]0,\,T[$.

(ii)  Assume that a point $q\in\Lambda_{\partial D}(p)$ is {\it
known} and there exists an open ball $U$ centered at $q$ such that
$U\cap\partial D$ is {\it known}. Extract the values of $\gamma$
and $\beta$ at $q$ from the data $u_f(x,t)$ given at all $x\in B$
and $t\in\,]0,\,T[$.

The aim of this paper is to give some solutions to (i) and (ii) by employing the {\it enclosure method}
as the {\it guiding principle}.
The enclosure method was introduced in \cite{I1,IF} and aims
at extracting a domain that {\it encloses} unknown {\it discontinuity}, such as, inclusion, cavity, crack,
obstacle, etc. embedded in a known background medium.  
It is based on, originally, the decaying and growing property of the {\it complex exponential
solutions} or the {\it complex geometrical optics solutions} of the governing equation 
of the used signal which propoagaes in the background medium.
The idea of the enclosure method used here goes back to \cite{I1}.
It is a {\it single measurement version} of the enclosure method.
Therein the governing equation is given by the Laplace equation and the idea yielded
an extraction formula of the support function for a polygonal cavity from a {\it single set} of the Cauchy data.
The formula has been tested numerically in \cite{IO} and the idea of this enclosure method
has been realized also for the inverse conductivity problem in two dimensions \cite{ICAL, IR},
inverse obstacle scattering problems at a fixed wave number in two dimensions
\cite{ISC, IH, ILOG, INS} and an extension of \cite{I1} to elastic body in two dimensions \cite{IITOU} and references therein.

Recently the single measurement version of the enclosure method
was developed also in \cite{I4} for the heat and wave equations in {\it one-space dimension}.
This paper opened the door to possibility of using several exponential solutions of the time dependent governing equation in 
the framework of that method.
Now we have already some results using this time domain single measurement version of the enclosure method in {\it three-space} dimensions
for the wave equation in \cite{IE0,IE3,IE4, IE5} and heat equation in \cite{IK2}.

This paper is closely related to some results in \cite{IE3, IE4}.
For (i) we have already known that from the data $u_f(x,t)$ given at all
$x\in B$ and $t\in\,]0,\,T[$ one can extract $\text{dist}\,(D,B)$ via the formula
$$\displaystyle
\lim_{\tau\longrightarrow\infty}\frac{1}{2\tau}\log\left\vert\int_B(w_f-v_f)dx\right\vert=
-\text{dist}\,(D,B),
\tag {1.2}
$$
where
$$\displaystyle
w_f(x,\tau)=\int_0^T e^{-\tau t}u_f(x,t)dt,\,\,x\in\Bbb R^3\setminus\overline D,\,\,\tau>0
$$
and $v_f\in H^1(\Bbb R^3)$ is the unique weak solution of the modified Helmholtz equation
$\displaystyle (\triangle-\tau^2)v+f=0$ in $\Bbb R^3$ which is given by
$$\displaystyle
v_f(x,\tau)=\frac{1}{4\pi}\int_{\Bbb R^3}\frac{e^{-\tau\vert x-y\vert}}{\vert x-y\vert}f(y)dy
=\frac{1}{4\pi}\int_B\frac{e^{-\tau\vert x-y\vert}}{\vert x-y\vert}dy.
$$
Note that the function $\tilde{v}(x,t)=e^{-\tau t}v_f(x,\tau)$
satisfies the inhomogeneous wave equation $\displaystyle
(\partial_t^2-\triangle)\tilde{v}+e^{-\tau t}f=0$ in $\Bbb
R^3\times\Bbb R$ and {\it decays} everwhere as $\tau\longrightarrow\infty$ unlike previous complex exponential solutions
of, e.g., the Laplace equation.

Therein it is assumed that
$\gamma(x)\le 1-C$ a.e. $x\in\partial D$ or $\gamma(x)\ge 1+C$
a.e. $x\in\partial D$ for a positive constant $C$ and the
reasonable constraint on $T$:
$$\displaystyle
T>2\text{dist}\,(D,B).
\tag {1.3}
$$ See Theorem 1.2 in \cite{IE3} for the detail.
Since $\text{dist}\,(D,B)=d_{\partial D}(p)-\eta$, the formula above yields $d_{\partial D}(p)$ which
is the radius of the largest sphere whose exterior {\it encloses} $D$.

As a corollary of (1.2), one gets a criterion whether given direction $\omega\in S^2$
the point $p+d_{\partial D}(p)\omega$ belongs to $\Lambda_{\partial D}(p)$ by using the back-scattering data $u_f$ on $B_{\eta-s}(p+s\omega)\times\,]0,\,T[$
for $f=\chi_{B_{\eta-s}(p+s\omega)}$ with a fixed $s\in\,]0,\,\eta[$, that is,
$p+d_{\partial D}(p)\omega\in\Lambda_{\partial D}(p)$ if and only if the quantity $d_{\partial D}(p+s\omega)$ computed by using the data
$u_f$ via the formula (1.2) with $f$ above coincides with $d_{\partial D}(p)-s$ (see Proposition 5.1 in \cite{IE4} for this type of a statement in the interior problem).

How about the shape of $D$?
This is the one of two questions considered in this paper.
Before describing a first result concerning with the question we give some remarks.

$\bullet$  If $q\in\Lambda_{\partial D}(p)$, then $q\in \partial B_{d_{\partial D}(p)}(p)$
and the two tangent planes at $q$ of $\partial D$ and $\partial B_{d_{\partial D}(p)}(p)$ coincide.

$\bullet$  We denote by $S_q(\partial D)$ and $S_q(\partial B_{d_{\partial D}(p)}(p))$ the
{\it shape operators} (or the {\it Weingarten maps})
at $q$ {\it with respect to $\nu_q$} which is the unit outward normal on $\partial D$ and inward on $\partial B_{d_{\partial D}(p)}(p)$.
Those are symmetric linear operators on the common tangent space at $q$ of $\partial D$ and $\partial B_{d_{\partial D}(p)}(p)$.

$\bullet$  $S_q(\partial B_{d_{\partial D}(p)}(p))-S_q(\partial D)\ge 0$ as the quadratic form
on the same tangent space at $q$ since $q$ attains $\min_{x\in\partial D}\vert x-p\vert$.

Now we can describe the following result which is the core of an answer to the question in the case when $\gamma\equiv 0$.

\proclaim{\noindent Theorem 1.1.}  Let $\gamma\equiv 0$.
Assume that $\partial D$ is $C^3$ and $\beta\in C^2(\partial D)$; $\Lambda_{\partial D}(p)$ is finite and satisfies
$$\displaystyle
\text{det}\,(S_q(\partial B_{d_{\partial D}(p)}(p))-S_q(\partial D))>0,\,\,\forall q\in\Lambda_{\partial D}(p).
$$
If $T$ satisfies (1.3),
then we have
$$\begin{array}{c}
\displaystyle
\lim_{\tau\longrightarrow\infty}
\tau^4e^{2\tau\text{dist}\,(D,B)}
\int_B(w_f-v_f)dx
=\frac{\pi}{2}
\left(\frac{\eta}{d_{\partial D}(p)}\right)^2A_{\partial D}(p),
\end{array}
\tag {1.4}
$$
where
$$\begin{array}{c}
\displaystyle
A_{\partial D}(p)
=\sum_{q\in\Lambda_{\partial D}(p)}
\frac{1}{\sqrt{\text{det}\,(S_q(\partial B_{d_{\partial D}(p)}(p))-S_q(\partial D))}}.
\end{array}
$$

\endproclaim

{\bf Remark 1.1.}
Theorem 1.1 tells us that
formula (1.4) is {\it invariant} with respect to the zeroth-order perturbation $\partial u/\partial\nu-\beta(x) u=0$
of the Neumann boundary condition $\partial u/\partial\nu=0$ on $\partial D$.
It seems that the proof of Theorem 1.1 cannot cover the case when $\gamma\not\equiv 0$.
The study for this case belongs to our next project.  See also \cite{M} for some results
using the {\it scattering amplitude} in the Lax-Phillips scattering theory when $\gamma\not\equiv 0$ and $\beta\equiv 0$.

{\bf Remark 1.2.}
Let $k_1(q)\le k_2(q)$ denote the eigenvalues of $S_q(\partial D)$.  They are called the {\it principle
curvatures} of $\partial D$ at $q$ {\it  with respect to} $\nu_q$.
Since $S_q(\partial B_{d_{\partial D}(p)}(p))=(1/d_{\partial D}(p))I$, we have
$$\displaystyle
\text{det}\,(S_q(\partial B_{d_{\partial D}(p)}(p))-S_q(\partial D))
=(\lambda-k_1(q))(\lambda-k_2(q)),
$$
where $\lambda=1/d_{\partial D}(p)$.  Recall the Gauss curvature $K_{\partial D}(q)$ of $\partial D$
at $q$ and mean curvature $H_{\partial D}(q)$ with respect to $\nu_q$ are given by
$$\displaystyle
K_{\partial D}(q)=k_1(q)k_2(q),\,\,H_{\partial D}(q)=(k_1(q)+k_2(q))/2.
$$
This yields another expression
$$\displaystyle
\text{det}\,(S_q(\partial B_{d_{\partial D}(p)}(p))-S_q(\partial D))
=\lambda^2-2H_{\partial D}(q)\lambda+K_{\partial D}(q),
\tag {1.5}
$$
where $\lambda$ the same as above.

As a corollary of Theorem 1.1 we have the following result which
enables us to extract information about the shape of $\partial D$
at a known $q\in \Lambda_{\partial D}(p)$ from two back-scattering data
corresponding to suitably chosen {\it two} initial data.

\proclaim{\noindent Corollary 1.1.}
Let $\gamma\equiv 0$. Assume that $\partial D$ is $C^3$ and $\beta\in C^2(\partial D)$.
Let $p\in\Bbb R^3\setminus\overline D$ and assume that $q\in\Lambda_{\partial D}(p)$ is known.
Let $B_1$ and $B_2$ denote two open balls cetred at $p-s_j\nu_q$, $j=1,2$, respectively with
$0<s_1<s_2<d_{\partial D}(p)$ and satisfy $\overline B_1\cup\overline B_2\subset\Bbb R^3\setminus\overline D$.

Then, one can extract $K_{\partial D}(q)$ and $H_{\partial D}(q)$
from $u_f$ on $B_j\times\,]0,\,T[$ with
$f=\chi_{B_j}$ for $j=1,2$ provided $T$ satisfies $T>2\max_{j=1,2}\text{dist}\,(D,B_j)$.
\endproclaim

Note that the centre of $B_2$ lies on the segment connecting the centre of $B_1$ with $q$;
$d_{\partial D}(p)=\vert p-q\vert$ and $\nu_q=(p-q)/\vert q-p\vert$.

The points are

$\bullet$ $\Lambda_{\partial D}(p-s_j\nu_q)=\{q\}$ and $\text{det}\,(S_q(\partial B_{d_{\partial D}(p-s_j\nu_q)}(p-s_j\nu_q))-S_q(\partial D))>0$;

$\bullet$ $d_{\partial D}(p-s_j\nu_q)=d_{\partial D}(p)-s_j$.

These enable us to apply Theorem 1.1 to the case when
$B=B_j$ and $f=\chi_{B_j}$ with $j=1,2$. Then with the help of (1.5), from (1.4) with
$j=1,2$ we have the following $2\times 2$ liner system for two
unknowns $K_{\partial D}(q)$ and $H_{\partial D}(q)$ via formula (1.4):
$$\begin{array}{c}
\displaystyle
Q(s)=(d_{\partial D}(p)-s)^{-2}-2H_{\partial D}(q)(d_{\partial D}(p)-s)^{-1}+K_{\partial D}(q),\,\,s=s_j,
\end{array}
$$
where $Q(s)$, $s=s_j$ are known quantities coming from (1.4).

By solving this system we obtain both the Gauss and mean curvatures of $\partial D$ at $q$.
Note that this idea goes back to the proof of Theorem 5.1 in \cite{IE4} where the interior problem
for the case when $\gamma=\beta=0$ has been considered.

By the way how about the question (ii)?  This is a new question in the framework of the enclosure method
and the complete answer remains open.
Clearly (1.3) does not help us since it does not contain any information about the coefficient $\beta$.
One possible direction is to clarify the remainder term of (1.4) as $\tau\longrightarrow\infty$, that is,
$$\displaystyle
\tau^4e^{2\tau\text{dist}\,(D,B)}
\int_B(w_f-v_f)dx
-\frac{\pi}{2}
\left(\frac{\eta}{d_{\partial D}(p)}\right)^2A_{\partial D}(p).
$$

The following theorem is closely related to this question and the main result of this paper.

\proclaim{\noindent Theorem 1.2.}  Let $\gamma\equiv 0$.
Assume that $\partial D$ is $C^5$ and $\beta\in C^2(\partial D)$;
$\Lambda_{\partial D}(p)$ is finite and satisfies
$$\displaystyle
\text{det}\,(S_q(\partial B_{d_{\partial D}(p)}(p))-S_q(\partial D))>0,\,\,\forall q\in\Lambda_{\partial D}(p).
$$
For each $q\in \Lambda_{\partial D}(p)$ let $\mbox{\boldmath
$e$}_j$, $j=1,2$ be orthogonal basis of the tangent space at $q$
of $\partial D$ with $\mbox{\boldmath $e$}_1\times\mbox{\boldmath
$e$}_2=\nu_q$. Choose an open ball $U$ centred at $q$ with radius
$r_q$ in such a way that: there exist a $h\in C^5_0(\Bbb R^2)$
with $h(0,0)=0$ and $\nabla h(0,0)=0$ such that $U\cap\partial
D=\{q+\sigma_1\mbox{\boldmath $e$}_1+\sigma_2\mbox{\boldmath
$e$}_2+h(\sigma_1,\sigma_2)\nu_q\,\vert\,
\sigma_1^2+\sigma_2^2+h(\sigma_1,\sigma_2)^2<r_q^2\}$.

If $T$ satisfies (1.3),
then we have
$$\begin{array}{c}
\displaystyle
\lim_{\tau\longrightarrow\infty}
\tau^5\left\{e^{2\tau\text{dist}\,(D,B)}
\int_B(w_f-v_f)dx-\frac{1}{\tau^4}\frac{\pi}{2}
\left(\frac{\eta}{d_{\partial D}(p)}\right)^2A_{\partial D}(p)\right\}\\
\\
\displaystyle
=-\frac{\pi\eta}{d_{\partial D}(p)^2}A_{\partial D}(p)
+\frac{\pi}{2}\eta^2B_{\partial D}(p),
\end{array}
\tag {1.6}
$$
where
$$\begin{array}{c}
\displaystyle
B_{\partial D}(p)
=\sum_{q\in\Lambda_{\partial D}(p)}
\frac{C_{\partial D}(q)}
{\displaystyle
\sqrt{\text{det}\,(S_q(\partial B_{d_{\partial D}(p)}(p))-S_q(\partial D))}},
\end{array}
$$
$$\begin{array}{c}
\displaystyle
C_{\partial D}(q)
=-\frac{1}{d_{\partial D}(p)^3}+\frac{11-12d_{\partial D}(p)H_{\partial D}(q)}
{\displaystyle
8d_{\partial D}(p)^5
\text{det}\,(S_q(\partial B_{d_{\partial D}(p)}(p))-S_q(\partial D))}\\
\\
\displaystyle
-\frac{1}{4d_{\partial D}(p)^2}
h_{\sigma_p\sigma_q\sigma_r}(0)h_{\sigma_s\sigma_t\sigma_u}(0)
\left(\frac{1}{4}B_{ps}B_{qr}B_{tu}+\frac{1}{6}B_{ps}B_{qt}B_{ru}\right)\\
\\
\displaystyle
+\frac{1}{16d_{\partial D}(p)^2}
h_{\sigma_p\sigma_q\sigma_r\sigma_s}(0)B_{pr}B_{qs}
-\frac{\beta(q)}{d_{\partial D}(p)^2}
\end{array}
$$
and
$$\displaystyle
B=(B_{pq})
=-\left(\frac{1}{d_{\partial D}(p)}I_2-\nabla^2h(0)\right)^{-1}.
$$

\endproclaim

From Theorem 1.2 we have immediately the following corollary.

\proclaim{\noindent Corollary 1.2.} Let $\gamma\equiv 0$.
Assume that $\partial D$ is $C^5$
and $\beta\in C^2(\partial D)$; $q\in\Lambda_{\partial D}(p)$ is
known; there exists an open ball $U$ centred at $q$ with radius
$r_q$ and orthonormal basis $\mbox{\boldmath $e$}_1$ and
$\mbox{\boldmath $e$}_2$ of the tangent space at $q$ of $\partial
D$ with $\mbox{\boldmath $e$}_1\times\mbox{\boldmath $e$}_2=\nu_q$
such that $U\cap\partial D=\{q+\sigma_1\mbox{\boldmath
$e$}_1+\sigma_2\mbox{\boldmath
$e$}_2+h(\sigma_1,\sigma_2)\nu_q\,\vert\,
\sigma_1^2+\sigma_2^2+h(\sigma_1,\sigma_2)^2<r_q^2\}$ with a $h\in
C^5_0(\Bbb R^3)$ satisfying $h(0,0)=0$, $\nabla h(0,0)=0$ and all
the second, third and fourth order derivatives of $h$ at
$\sigma=(0,0)$ are known.
Let $0<s<d_{\partial D}(p)$.
Let $B'$ denote the open ball centred at $p-s\nu_q$ and satisfy
$\overline B'\cap\overline D=\emptyset$.
Then, one can extract $\beta(q)$ from $u_f$ on $B'$ for
$f=\chi_{B'}$ provided $T$ satisfies $T>2\text{dist}\,(D,B')$.
\endproclaim

This corollary says, in short, assume that we know in advance, a
point $q\in \Lambda_{\partial D}(p)$ and thus $d_{\partial
D}(p)=\vert p-q\vert$ and $\nu_q=(p-q)/\vert p-q\vert$, too.  Thus
the tangent plane $(x-q)\cdot\nu_q=0$ at $q$ of $\partial D$ is
known. Moreover, assume that: we know that $\partial D$ near $q$
is given by making a {\it rotation} around the normal at $q$ of a
{\it graph} of a function $h$ defined on the tangent plane and
that, in an appropriate orthogonal coordinates on the tangent
plane, say $\sigma=(\sigma_1,\sigma_2)$, the Taylor expansion of
the function at $\sigma=0$ has the form
$h(\sigma_1,\sigma_2)=\sum_{2\le \vert\alpha\vert\le
4}h_{\alpha}\sigma^{\alpha}+\cdots$ with {\it known} coefficients
$h_{\alpha}$ for $2\le\vert\alpha\vert\le 4$. Then, produce the
wave $u_f$ with $f=\chi_{B'}$ for a small $s>0$
and measure the wave on $B'$. Since
$\Lambda_{\partial D}(p-s\nu_q)=\{q\}$, computing (1.6) in which
$B$ is replaced with $B'$ and $p$ with $p-s\nu_q$, we obtain
$C_{\partial D}(q)$ which is a linear equation with unknown $\beta(q)$ and thus
solving this, one obtains $\beta(q)$.

By the way, it is also interesting to consider the following
inverse problem.

{\bf Problem II.}
Let $B$ and $f$ be same as those of Problem I.

Let $B_R(p)$ denote the open ball centred at $p$ with radius $R$
and satisfy $B_R(p)\subset\Bbb R^3\setminus\overline D$ and $B\subset B_R(p)$.

(i)'  Assume that $D$ is {\it unknown} and both $\gamma$ and $\beta$ {\it unknown}.  Extract information
about the location and shape of $D$ from the data $u_f(x,t)$ given at all
$x\in\partial B_D(p)$ and $t\in\,]0,\,T[$.

(ii)'  Assume that a point $q\in\Lambda_{\partial D}(p)$ is {\it
known} and there exists an open ball $U$ centered at $q$ such that
$U\cap\partial D$ is {\it known}. Extract the values of $\gamma$
and $\beta$ at $q$ from the data $u_f(x,t)$ given at all $x\in\partial B_R(p)$
and $t\in\,]0,\,T[$.

The difference between Problem I and II is the place where the wave is observed.
The latter case is clearly desirable since the dimension of the observation place is lower
than that of the former case.  The sphere $\partial B_R(p)$ is a model of the place where many receivers
are placed.

For this problem we show that the enclosure method does not make an issue of this difference at all.
The point is the following asymptotic formula between two data.

\proclaim{\noindent Proposition 1.1.}
Let $B$ and $f$ be same as those of Problem I.
Let $B_R(p)$ denote the open ball centred at $p$ with radius $R$
and satisfy $B_R(p)\subset\Bbb R^3\setminus\overline D$ and $B\subset B_R(p)$.

We have, as $\tau\longrightarrow\infty$
$$\begin{array}{c}
\displaystyle
\int_{\partial B_R(p)}(w_f-v_f)dS
\\
\\
\displaystyle
=
\frac{R}{\eta}\tau e^{\tau(R-\eta)}(1+O(\tau^{-1}e^{-\tau\eta}))
\int_B(w_f-v_f)dx+O(\tau e^{-\tau(T-(R-\eta))}).
\end{array}
\tag {1.7}
$$

\endproclaim

Note that $w_f$ on $\partial B_R(p)$ can be computed from the trace of $u_f(\,\cdot\,,t)$ onto $\partial B_R(p)$
given at almost all $t\in\,]0,\,T[$ via the formula
$$\displaystyle
w_f=\int_0^T e^{-\tau t}u_f(\,\cdot\,,t)\vert_{\partial B_R(p)}dt\in H^{1/2}(\partial B_R(p)).
$$
We can say that throughout formula (1.7) all the results mentioned
above are transplanted in this case. For the detailed description
of the transplanted results see Section 4.

A brief outline of the rest of the paper is as follows.
In section 2 we derive an asymptotic representation formula of the {\it indicator function}
$$\displaystyle
\tau\longmapsto \int_B(w_f-v_f)dx
$$
as $\tau\longrightarrow\infty$.  It clarifies the {\it principle term} in terms of integrals
of $v_f$ over $D$ and $\partial D$.  The key for the proof is an argument developed in \cite{LP} which we call the Lax-Phillips
{\it reflection argument}.  Previously, we applied the argument
for the case when $\gamma\equiv 0$ and $\beta\equiv 0$ in \cite{IE4}.
In this paper, we found that some {\it modification} of the argument still works in our problem setting.

Having the formula established in section 2, we prove theorems 1.1
and 1.2 in Section 3. They are an application of the Laplace
method, however, we need the second term of the asymptotic
expansion of a double integral of Laplace type. As can be seen in
\cite{BH}, the coefficient is quite complicated. To make the
relation between the obstacle and the coefficient as concise as
possible, we did some additional calculation and they are
summarized as Lemma 3.1. Although the calculation is simple, it is
tedious. Thus, we put the proof of Lemma 3.1 in the Appendix.
In Section 4 we give a proof of Proposition 1.1 and its implications to Problem II.

Befor closing Introduction, we give a remark. In the Lax-Phillips
scattering theory, because of the difference of the purpose from
us, the wave is observed far a way from the obstacle and thus,
{\it infinite} observation time is needed. In the case when
$\gamma\not\equiv 0$ and $\beta\equiv 0$, there is a classical
result due to Majda \cite{M} using the high frequency asymptotics
of the scattering amplitude which is the Fourier transform of the
{\it scattering kernel} in the Lax-Phillips scattering theory. The
method of the proof used therein is based on the idea of
geometrical optics and thus, different from us. And our result
yields not only the Gauss curvature but also the mean curvature.
This is also another difference from his result.

\section{Extracting the principle term of the indicator function}

In this section we set $c=c(x,\tau)=\gamma(x)\tau+\beta(x)$.
Let $f\in L^2(\Bbb R^3)$ satisfy $\text{supp}\,f\subset B$.

We give a proof of the following asymptotic formula.

\proclaim{\noindent Proposition 2.1.}  Let $\gamma\equiv 0$.
Assume that $\partial D$ is $C^3$ and $\beta\in C^2(\partial D)$.
We have
$$\displaystyle
\int_Bf(w_f-v_f)dx
=2J(\tau)(1+O(\tau^{-1/2}))+O(\tau^{-1}e^{-\tau T}),
\tag {2.1}
$$
where
$$\displaystyle
J(\tau)
=\int_D(\vert\nabla v_f\vert^2+\tau^2\vert v_f\vert^2)dx
-\int_{\partial D}c\vert v_f\vert^2 dS.
$$

\endproclaim

Let $\epsilon_f^0\in H^1(\Bbb R^3\setminus\overline D)$ solve
$$\begin{array}{c}
\displaystyle
(\triangle-\tau^2)\epsilon_f^0=0\,\,\text{in}\,\Bbb R^3\setminus\overline D,\\
\\
\displaystyle
\frac{\partial\epsilon_f^0}{\partial\nu}
-c\epsilon_f^0=-\left(\frac{\partial v_f}{\partial\nu}
-cv_f\right)\,\,\text{on}\,\partial D.
\end{array}
\tag {2.2}
$$
and define
$$\displaystyle
E(\tau)
=\int_{\Bbb R^3\setminus\overline D}
(\vert\nabla\epsilon_f^0\vert^2+\tau^2\vert\epsilon_f^0\vert^2)dx
+\int_{\partial D}c\vert\epsilon_f^0\vert^2 dS.
$$

Clearly Proposition 2.1 is a consequence of two asymptotic formulae (2.3) and (2.8)
described in the following two lemmas.
The core of this section is the proof of (2.8) which forms a subsection
independently and can be considered as the one of essential parts of this paper.

\proclaim{\noindent Lemma 2.1.}
As $\tau\longrightarrow\infty$ we have
$$\displaystyle
\int_Bf(w_f-v_f)dx
=J(\tau)+E(\tau)+O(\tau^2 e^{-\tau (T+\text{dist}\,(D,B))})+O(\tau^{-1}e^{-\tau T}).
\tag {2.3}
$$

\endproclaim

{\it\noindent Proof.}
Integration by parts yields
$$\begin{array}{c}
\displaystyle
\int_Bf(w_f-v_f)dx
=\int_{\partial D}\left(\frac{\partial v}{\partial\nu}-cv\right)wdS
-e^{-\tau T}
\left(\int_{\Bbb R^3\setminus\overline D}Fvdx+\int_{\partial D}GvdS\right),
\end{array}
\tag {2.4}
$$
where
$$\displaystyle
F=F(x,\tau)=\partial_tu(x,T)+\tau u(x,T)
$$
and
$$\displaystyle
G=G(x)=\gamma(x)u(x,T).
$$
Note that the definition of the weak solution taken from \cite{DL} ensures $\Vert F\Vert_{L^2(\Bbb R^3\setminus\overline D)}=O(\tau)$ and
$\Vert G\Vert_{L^2(\partial D)}<\infty$ (see \cite{IE3}).
Since $v_f$ has the bounds
$\displaystyle
\Vert v_f\Vert_{L^2(\Bbb R^3)}=O(\tau^{-2})$
and $\Vert v\Vert_{L^{\infty}(\partial D)}=O(e^{-\tau\text{dist}\,(D,B)})$, we see that the second term of (2.4) has the bound
$O(\tau^{-1}e^{-\tau T})+O(e^{-\tau(T+\text{dist}\,(D,B))})$.

One can write $w_f=v_f+\epsilon_f^0+Z$, where $Z$ solves
$$\begin{array}{c}
\displaystyle
(\triangle-\tau^2)Z=e^{-\tau T}F\,\,\text{in}\,\Bbb R^3\setminus\overline D,\\
\\
\displaystyle
\frac{\partial Z}{\partial\nu}-cZ=e^{-\tau T}G\,\,\text{on}\,\partial D.
\end{array}
$$
Using the same argument as done in the proof of Lemma 2.1 in \cite{IE3}, one has
$$\begin{array}{c}
\displaystyle
\int_{\Bbb R^3\setminus\overline D}(\vert\nabla Z\vert^2+\tau^2\vert Z\vert^2)dx
=O(\tau^2e^{-2\tau T}).
\end{array}
\tag {2.5}
$$
A combination of this, the trace theorem onto $\partial D$ and the
trivial estimate $\Vert\partial
v_f/\partial\nu-cv_f\Vert_{L^{\infty}(\partial D)}=O(\tau
e^{-\tau\text{dist}\,(D,B)})$, we obtain
$$\displaystyle
\int_{\partial D}\left(\frac{\partial v_f}{\partial\nu}-cv_f\right)ZdS=O(\tau^2 e^{-\tau (T+\text{dist}\,(D,B))}).
$$
It is easy to see that integration by parts yields
$$\displaystyle
J(\tau)
=\int_{\partial D}
\left(\frac{\partial v_f}{\partial\nu}-cv_f\right)v_fdS
\tag {2.6}
$$
and
$$\displaystyle
E(\tau)
=\int_{\partial D}\left(\frac{\partial v_f}{\partial\nu}
-cv\right)\epsilon_f^0dS.
\tag {2.7}
$$
Now (2.3) is clear.

\noindent
$\Box$

\proclaim{\noindent Lemma 2.2.}
Let $\gamma\equiv 0$ and $\beta\in C^2(\partial D)$.  There exists a $\tau_0>0$ such that, for all $\tau\ge\tau_0$ we have $J(\tau)>0$
and as $\tau\longrightarrow\infty$,
$$\displaystyle
E(\tau)=J(\tau)(1+O(\tau^{-1/2})).
\tag {2.8}
$$

\endproclaim

\subsection{Proof of Lemma 2.2}

The proof is a modification of the Lax-Phillips reflection argument \cite{LP}(see also \cite{IE4} for the case when $c\equiv 0$).
In this subsection we always assume that $\partial D$ is $C^3$ and $c\in C^2(\partial D)$.
And for simplicity we denote $v_f$ and $\epsilon_0$ by $v$ and $\epsilon_0$, respectively. 

\subsubsection{A representation formula of $E(\tau)-J(\tau)$ via a reflection}

There exists a positive constant
$\delta_0$ such that: given $x\in \Bbb
R^3\setminus D$/$x\in\overline D$ with $d_{\partial D}(x)<2\delta_0$ there
exists a unique $q=q(x)$ be the boundary point on $\partial D$
such that $x=q\pm d_{\partial D}(x)\nu_q$ (\cite{GT}).
Both $d_{\partial D}(x)$ and $q(x)$ are $C^2$
for $x\in\Bbb R^3\setminus D$/$x\in D$ with $d_{\partial D}(x)<2\delta_0$
(Lemma 1 of Appendix in \cite{GT}).
Note that
$\nu_q$ is the unit outer normal to $\partial D$ at $q$.

For $x$ with $d_{\partial D}(x)<2\delta_0$ define
$$\displaystyle
x^r=2q(x)-x.
$$
Let $0<\delta<\delta_0$.
Let $\phi=\phi_{\delta}$ be a smooth cut-off
function, $0\le\phi(x)\le 1$, and such that: $\phi(x)=1$ if
$d_{\partial D}(x)<\delta$ and $\phi(x)=0$ if $d_{\partial
D}(x)>2\delta$;
$\vert \nabla\phi(x)\vert\le C\delta^{-1}$;
$\vert\nabla^2\phi(x)\vert\le C\delta^{-2}$.

Define
$$\displaystyle
v^r(x)=\phi(x)v(x^r),\,x\in\Bbb R^3\setminus\overline D.
$$
Let $\eta=\eta(x)$ be a $C^2$-function of $x\in\Bbb R^3\setminus D$
with $d_{\partial D}(x)<2\delta_0$.

\proclaim{\noindent Proposition 2.2.} Assume that $\partial D$ is
$C^3$ and that $\gamma\in C^2(\partial D)$ and $\beta\in
C^2(\partial D)$. Let $\eta$ satisfy $\eta=1$ on $\partial D$ and
$$
\displaystyle
\frac{\partial\eta}{\partial\nu}=2(\gamma(x)\tau+\beta(x))\,\,\text{on}\,\partial D.
$$
Then, we have
$$\begin{array}{c}
\displaystyle
E(\tau)-J(\tau)=
\int_{\Bbb R^3\setminus\overline D}\epsilon_0(\triangle-\tau^2)(\eta v^r)dx.
\end{array}
\tag {2.9}
$$

\endproclaim

{\it\noindent Proof.}
Integration by parts (or the weak formulation of (2.2)) yields
$$\begin{array}{c}
\displaystyle
\int_{\Bbb R^3\setminus\overline D}\epsilon_0(\triangle-\tau^2)(\eta v^r)dx
=\int_{\partial D}
\left\{\frac{\partial\epsilon_0}{\partial\nu}(\eta v^r)
-\epsilon_0\frac{\partial}{\partial\nu}(\eta v^r)\right\}dS.
\end{array}
\tag {2.10}
$$
Since
$$\begin{array}{c}
\displaystyle
\frac{\partial v^r}{\partial\nu}=-\frac{\partial v}{\partial\nu}\,\,,
v^r=v\,\,\text{on}\,\partial D,
\end{array}
$$
we have
$$\displaystyle
\frac{\partial}{\partial\nu}(\eta v^r)
=\frac{\partial\eta}{\partial\nu}v-\eta\frac{\partial v}{\partial\nu}\,\,\text{on}\,\partial D.
$$
Substituting this together with
$$\displaystyle
\frac{\partial\epsilon_0}{\partial\nu}
=c\epsilon_0-\left(\frac{\partial v}{\partial\nu}-cv\right)\,\,\text{on}\,\partial D
$$
into (2.10), we obtain
$$\begin{array}{c}
\displaystyle
\int_{\Bbb R^3\setminus\overline D}\epsilon_0(\triangle-\tau^2)(\eta v^r)dx
=
\int_{\partial D}
\left\{
\left(c\epsilon_0-\left(\frac{\partial v}{\partial\nu}-cv\right)\right)
\eta v
-\epsilon_0
\left(\frac{\partial\eta}{\partial\nu}v-\eta\frac{\partial v}{\partial\nu}\right)
\right\}dS
\\
\\
\displaystyle
=\int_{\partial D}
\left\{\left(c\eta-\frac{\partial\eta}{\partial\nu}\right)v+
\eta\frac{\partial v}{\partial\nu}\right\}\epsilon_0dS
-\int_{\partial D}\left(\frac{\partial v}{\partial\nu}-cv\right)\eta vdS\\
\\
\displaystyle
=\int_{\partial D}
\left\{\frac{\partial v}{\partial\nu}-\frac{1}{\eta}
\left(\frac{\partial\eta}{\partial\nu}-c\eta\right)v\right\}\eta\epsilon_0dS
-\int_{\partial D}
\left(\frac{\partial v}{\partial\nu}-cv\right)\eta vdS.
\end{array}
$$
This together with (2.6), (2.7) yields (2.9).

\noindent
$\Box$

One possible choice of $\eta$ in (2.9) is
$$\displaystyle
\eta(x)
=1+2c(q(x),\tau)d_{\partial D}(x),\,\,x\in\Bbb R^3\setminus D,\,\,d_{\partial D}(x)<2\delta_0
\tag {2.11}
$$
since we have
$$\displaystyle
\nabla(d_{\partial D}(x))=\nu(q(x))
$$
and
$$\displaystyle
\nabla q^j(x)\cdot\nu(q(x))=0,\,j=1,2,3.
$$

In what follows we always make use of $\eta$ given by (2.11) and thus we have (2.9) for this $\eta$.

\subsubsection{The Lax-Phillips reflection argument}

We assume that $\gamma\equiv 0$ and $\beta\in C^2(\partial D)$.  Thus (2.11) becomes
$\eta(x)=\beta(q(x))d_{\partial D}(x)$ with $d_{\partial D}(x)<2\delta_0$.

In this subsubsection we give the following upper bound for the integral in the right-hand side on (2.9):
there exists a $\tau_0>0$ such that, for all $\tau\ge\tau_0$ $J(\tau)>0$ and $E(\tau)>0$
and as $\tau\longrightarrow\infty$
$$\displaystyle
\left\vert
\int_{\Bbb R^3\setminus\overline D}\epsilon_0(\triangle-\tau^2)(\eta v^r)dx
\right\vert
\le O(\tau^{-1/2})(E(\tau)J(\tau))^{1/2}.
\tag {2.12}
$$
If once we have this, from (2.9) we have
$$\displaystyle
\vert E(\tau)-J(\tau)\vert\le O(\tau^{-1/2}(E(\tau)J(\tau))^{1/2}.
\tag {2.13}
$$
Then, taking the square of the both sides on (2.13), we obtain
$$\displaystyle
E^2(\tau)\le (2+O(\tau^{-1}))E(\tau)J(\tau)
$$
and thus $E(\tau)\le (2+O(\tau^{-1}))J(\tau)$.
A combination of this and (2.13) yields
$$\displaystyle
\vert E(\tau)-J(\tau)\vert\le O(\tau^{-1/2})J(\tau).
$$
Therefore we obtain (2.8).  This completes the proof of Lemma 2.2.

Thus the point is the derivation of (2.13). Here we apply the
Lax-Phillips reflection argument developed in \cite{LP}(see also
Appendix A in \cite{IE4}) to our situation.  We focus only on the
essential part of the derivation.

The key point is  a differential identity for $(\triangle-\tau^2)(v^r)$
which is a consequence of (4.15) in \cite{LP}.
It takes the form
$$\begin{array}{c}
\displaystyle
(\triangle-\tau^2)(v^r)(x)
=\phi(x)\sum_{i,j}d_{\partial D}(x)a_{ij}(x)(\partial_i\partial_jv)(x^r)
\\
\\
\displaystyle
+\sum_{j}\left(\sum_{k}b_{jk}(x)\frac{\partial\phi}{\partial x_k}(x)+d_j(x)\phi(x)\right)\frac{\partial v}{\partial x_j}(x^r)
+(\triangle\phi)(x)v(x^r),
\end{array}
$$
where
$a_{ij}(x)$, $b_{j,k}(x)$ and $d_j(x)$ with $i,j,k=1,2,3$ are independent of $\phi$ and $v$; $a_{ij}(x)$ and
$b_{jk}(x)$ are $C^1$ and $d_j(x)$ is $C^0$ for $x\in\Bbb R^3\setminus D$ with $d_{\partial D}(x)<2\delta_0$.

A change of variables $x=y^r$ gives
$$\begin{array}{c}
\displaystyle
\int_{\Bbb R^3\setminus\overline D}\epsilon_0(x)\eta(x)(\triangle-\tau^2)v^r(x)dx\\
\\
\displaystyle
=\int_{D}\epsilon_0(y^r)\eta(y^r)\left\{\phi(y^r)\sum_{i,j}d_{\partial D}(y^r)a_{ij}(y^r)(\partial_i\partial_jv)(y)+\text{lower order terms}\right\}
J(y)dy,
\end{array}
\tag {2.14}
$$
where $J(y)$ denotes the Jacobian.  The point is the bound on the first term involving the second order derivatives of $v$.
Since $d_{\partial D}(x)=0$ on $\partial D$, applying the integration by parts to the term, one can rewrite
$$\begin{array}{c}
\displaystyle
\int_{D}\epsilon_0(y^r)\eta(y^r)\phi(y^r)\sum_{i,j}d_{\partial D}(y^r)a_{ij}(y^r)(\partial_i\partial_jv)(y)J(y)dy\\
\\
\displaystyle
=-\sum_{i,j}\int_{D}\partial_i\{\epsilon_0(y^r)\eta(y^r)\phi(y^r)d_{\partial D}(y)a_{ij}(y^r)J(y)\}\partial_jv(y)dy.
\end{array}
$$
Note that $d_{\partial D}(y^r)=d_{\partial D}(y)$.
This yields
$$\begin{array}{c}
\displaystyle
\left\vert\int_{D}\epsilon_0(y^r)\eta(y^r)\phi(y^r)\sum_{i,j}d_{\partial D}(y^r)a_{ij}(y^r)(\partial_i\partial_jv)(y)J(y)dy\right\vert\\
\\
\displaystyle
=\left\{O(\delta)\Vert\nabla\epsilon_0^r\Vert_{L^2(D_{\delta})}+O(1)\Vert\epsilon_0^r\Vert_{L^2(D_{\delta})})
\right\}\Vert\nabla v\Vert_{L^2(D)},
\end{array}
\tag {2.15}
$$
where $\epsilon_0^r(y)=\epsilon_0(y^r)$ and $D_{\delta}=\{y\in
D\,\vert\,d_{\partial D}(y)<2\delta\}$. Note that $O(\delta)$ is
coming from $d_{\partial D}(y)$. Just simply
estimating the integrals involving the lower order terms in the
right-hand side on (2.14) and combining the results with (2.15),
we obtain
$$\begin{array}{c}
\displaystyle
\left\vert\int_{\Bbb R^3\setminus\overline D}\epsilon_0(x)\eta(x)(\triangle-\tau^2)v^r(x)dx\right\vert=O(\delta)\Vert\nabla\epsilon_0^r\Vert_{L^2(D_{\delta})}\Vert\nabla v\Vert_{L^2(D)}\\
\\
\displaystyle
+O(\delta^{-1})\Vert \epsilon_0^r\Vert_{L^2(D_{\delta})}\Vert\nabla v\Vert_{L^2(D)}
+O(\delta^{-2})\Vert \epsilon_0^r\Vert_{L^2(D_{\delta})}\Vert v\Vert_{L^2(D)}.
\end{array}
\tag {2.16}
$$

It is easy to see that
$$
\displaystyle
\Vert \epsilon_0^r\Vert_{L^2(D_{\delta})}
\le C\Vert \epsilon_0\Vert_{L^2(\Bbb R^3\setminus\overline D)},\,\,
\Vert\nabla \epsilon_0^r\Vert_{L^2(D_{\delta})}\le
C\Vert\nabla \epsilon_0\Vert_{L^2(\Bbb R^3\setminus\overline D)}.
\tag {2.17}
$$
Moreover, using a trace theorem, one can easily obtain the following
inequalities: there exist constants $\tau_0>0$ and $C>0$ such
that, for all $\tau\ge\tau_0$
$$\displaystyle
\Vert\epsilon_0\Vert_{L^2(\Bbb R^3\setminus\overline D)}^2\le C\tau^{-2}E(\tau),
\tag {2.18}
$$
$$\displaystyle
\Vert\nabla\epsilon_0\Vert_{L^2(\Bbb R^3\setminus\overline D)}^2\le C E(\tau),
\tag {2.19}
$$
$$\displaystyle
\Vert v\Vert_{L^2(D)}^2\le C\tau^{-2}J(\tau),
\tag {2.20}
$$
$$\displaystyle
\Vert\nabla v\Vert_{L^2(D)}^2\le CJ(\tau).
\tag {2.21}
$$
Note that these include also $E(\tau)>0$ and $J(\tau)>0$ for all $\tau\ge\tau_0$ provided $\gamma\equiv 0$.
Applying (2.18) and (2.19) to (2.17), we obtain
$$\displaystyle
\Vert \epsilon_0^r\Vert_{L^2(D_{\delta})}\le C\tau^{-1}E(\tau)^{1/2},\,\,
\Vert\nabla \epsilon_0^r\Vert_{L^2(D_{\delta})}\le CE(\tau)^{1/2}.
\tag {2.22}
$$
Now applying (2.20), (2.21) and (2.22) to the right-hand side on (2.16),
we obtain
$$\displaystyle
\left\vert\int_{\Bbb R^3\setminus\overline D}\epsilon_0(x)\eta(x)(\triangle-\tau^2)v^r(x)dx\right\vert
\le\{O(\delta)+O(\delta^{-1}\tau^{-1})+O(\delta^{-2}\tau^{-2})\}(E(\tau)J(\tau))^{1/2}.
\tag {2.23}
$$

On the other hand, a direct computation yields
$$\begin{array}{c}
\displaystyle
(\triangle-\tau^2)(\eta v^r)
-\eta(\triangle-\tau^2)(v^r)\\
\\
\displaystyle
=(\triangle\eta)(x)\phi(x)v(x^r)
+2(\nabla\eta(x)\cdot\nabla\phi(x))v(x^r)\\
\\
\displaystyle
+2\phi(x)(2q'(x)-I_3)\nabla\eta(x)\cdot(\nabla v)(x^r).
\end{array}
$$
From this and the regularity of $q(x)$ and $d_{\partial D}(x)$ for $x\in\Bbb R^3\setminus D$ with $d_{\partial D}(x)<2\delta_0$ in \cite{GT}
one obtains that:
there exists a positive constant $C$ such that, for all $x\in\Bbb R^3\setminus D$
$$\begin{array}{c}
\displaystyle
\left\vert
(\triangle-\tau^2)(\eta v^r)(x)
-\eta(x)(\triangle-\tau^2)(v^r)(x)\right\vert\\
\\
\displaystyle
\le C\left\{\vert v(x^r)\vert(\vert\phi(x)\vert+
\vert\nabla\phi(x)\vert)
+\vert(\nabla v)(x^r)\vert\vert\phi(x)\vert\right\}.
\end{array}
\tag {2.24}
$$
Using the change of variables $x=y^r$,  we have also from (2.20), (2.21) and (2.22)
$$\begin{array}{c}
\displaystyle
\int_{\Bbb R^3\setminus\overline D}\vert\epsilon_0(x)\vert\left\{\vert v(x^r)\vert(\vert\phi(x)\vert+
\vert\nabla\phi(x)\vert)
+\vert(\nabla v)(x^r)\vert\vert\phi(x)\vert\right\}dx\\
\\
\displaystyle
\le C(\tau^{-2}(1+O(\delta^{-1}))+\tau^{-1})(E(\tau)J(\tau))^{1/2}.
\end{array}
$$
From this, (2.23) and (2.24) we obtain
$$\begin{array}{c}
\displaystyle
\left\vert
\int_{\Bbb R^3\setminus\overline D}\epsilon_0(\triangle-\tau^2)(\eta v^r)dx
\right\vert
\\
\\
\displaystyle
\le
\{O(\delta)+O(\delta^{-1}\tau^{-1})+O(\delta^{-2}\tau^{-2})
+C(\tau^{-2}(1+O(\delta^{-1}))+\tau^{-1})\}(E(\tau)J(\tau))^{1/2}.
\end{array}
$$
Now choosing $\delta=\tau^{-1/2}$, we obtain (2.12).

{\bf\noindent Remark 2.1.}  If $\gamma\not\equiv 0$, then the argument above does not work because of
the presence of the growing factor $\gamma(q(x))\tau$ as $\tau\longrightarrow\infty$.
Note also that $J(\tau)<0$ for $\tau>>1$ when $\inf_{x\in\partial D}(\gamma(x)-1)>0$ (see \cite{IE3}).

\section{Proof of Theorems 1.1 and 1.2}

In this section, we set $f=\chi_B$ and $v=v_f$.
By (2.1), it suffices to study the asymptotic behaviour of
$$\displaystyle
J(\tau)
=\int_{\partial D}
\left(\frac{\partial v}{\partial\nu}-cv\right)vdS\\
\\
=\int_{\partial D}\frac{\partial v}{\partial\nu}vdS-\int_{\partial D}c\vert v\vert^2dS.
\tag {3.1}
$$
We will see that the {\it second order term} not the top term contains information about $c=\beta(x)$ at $x\in\Lambda_{\partial D}(p)$.

By (3.26) in \cite{IE4}, we have, for a positive constant $C_1$,
$$\begin{array}{c}
\displaystyle
\int_{\partial D}\frac{\partial v}{\partial\nu}vdS
=\frac{1}{4\tau^3}
\left(\eta-\frac{1}{\tau}\right)^2
\int_{\partial D}
\frac{e^{-2\tau d_B(x)}}{\vert x-p\vert^2}
\left(1+\frac{1}{\tau\vert x-p\vert}\right)
\frac{p-x}{\vert p-x\vert}\cdot\nu_xdS\\
\\
\displaystyle
+O(\tau^{-1}e^{-2\tau\text{dist}\,(D,B)(1+C_1)}).
\end{array}
\tag {3.2}
$$
Note that $\nu_x$ is outward to $D$.

On the other hand, by (3.22) and (3.24) in \cite{IE4}, we have, for a positive constant $C_2$,
$$\displaystyle
\int_{\partial D}c\vert v\vert^2 dS
=\frac{1}{4\tau^4}
\left(\eta-\frac{1}{\tau}\right)^2
\int_{\partial D}
\frac{\beta(x)e^{-2\tau d_B(x)}}
{\vert x-p\vert^2} dS_x
+O(\tau^{-5}e^{-2\tau\text{dist}\,(D,B)(1+C_2)}).
\tag {3.3}
$$

Substituting (3.2) and (3.3) into (3.1), we obtain
$$\displaystyle
J(\tau)
=\frac{1}{4\tau^3}
\left(\eta-\frac{1}{\tau}\right)^2\left(I_1(\tau)+\frac{1}{\tau}I_2(\tau)\right)
+O(\tau^{-1}e^{-2\tau\text{dist}\,(D,B)(1+C_1)}),
\tag {3.4}
$$
where
$$\displaystyle
I_1(\tau)
=\int_{\partial D}\frac{1}{\vert x-p\vert^2}\frac{(p-x)}{\vert x-p\vert}\cdot\nu_x\,
e^{-2\tau d_B(x)}
dS_x
$$
and
$$\displaystyle
I_2(\tau)
=\int_{\partial D}
\left(\frac{1}{\vert x-p\vert^3}\frac{(p-x)}{\vert x-p\vert}\cdot\nu_x
-\frac{\beta (x)}{\vert x-p\vert^2}\right)e^{-2\tau d_B(x)}
dS_x.
$$

A combination of (2.1) and (3.4) gives
$$\begin{array}{c}
\displaystyle
\int_B(w_f-v_f)dx\\
\\
\displaystyle
=\left\{\frac{1}{2\tau^3}
\left(\eta-\frac{1}{\tau}\right)^2
\left(I_1(\tau)+\frac{1}{\tau}I_2(\tau)\right)
+O(\tau^{-1}e^{-2\tau\text{dist}\,(D,B)(1+C_1)})\right\}(1+O(\tau^{-1/2}))\\
\\
\displaystyle
+O(\tau^{-1}e^{-\tau T}).
\end{array}
\tag {3.5}
$$

{\bf Remark 3.1.}
Note that, to make the derivation of (3.4) self-contained, an alternative direct proof is also given in Appendix.

Now we study the asymptotic behaviour of $I_1(\tau)$ and $I_2(\tau)$ as $\tau\longrightarrow\infty$.

For each $q\in \Lambda_{\partial D}(p)$ let $\mbox{\boldmath
$e$}_j$, $j=1,2$ be orthogonal basis of the tangent space at $q$
of $\partial D$ with $\mbox{\boldmath $e$}_1\times\mbox{\boldmath
$e$}_2=\nu_q$. Choose an open ball $U$ centred at $q$ with radius
$r_q$ in such a way that: there exist a $h\in C^5_0(\Bbb R^2)$
with $h(0,0)=0$ and $\nabla h(0,0)=0$ such that $U\cap\partial
D=\{q+\sigma_1\mbox{\boldmath $e$}_1+\sigma_2\mbox{\boldmath
$e$}_2+h(\sigma_1,\sigma_2)\nu_q\,\vert\,
\sigma_1^2+\sigma_2^2+h(\sigma_1,\sigma_2)^2<r_q^2\}$. In
$U\cap\partial D$, we have
$$\displaystyle
x-p=\sigma\mbox{\boldmath $e$}_1+\sigma_2\mbox{\boldmath $e$}_2+(h(\sigma)-d_{\partial D}(p))\nu_q
$$
and
$$\displaystyle
\nu_x\,dS_x
=(-h_{\sigma_1}(\sigma)\mbox{\boldmath $e$}_1-h_{\sigma_2}(\sigma)\mbox{\boldmath $e$}_2+\nu_q)\,d\sigma.
$$
Thus we have
$$\displaystyle
\frac{p-x}{\vert x-p\vert}\cdot\nu_xdS_x
=\Psi_q(\sigma)^{-1}(\nabla h(\sigma)\cdot\sigma+d_{\partial D}(p)-h(\sigma))d\sigma,
$$
where
$$\displaystyle
\Psi_q(\sigma)=\vert x-p\vert
=\sqrt{\vert\sigma\vert^2+(d_{\partial D}(p)-h(\sigma))^2}.
$$
Define
$$\displaystyle
\phi(\sigma)=\Psi_q(0)-\Psi_q(\sigma).
$$
We have
$$\displaystyle
\phi_{\sigma_i\sigma_j}(0)
=-\left(\frac{1}{d_{\partial D}(p)}I_2
-\nabla^2h(0)\right)
$$
and thus
$$\displaystyle
\text{det}\,(\phi_{\sigma_i\sigma_j}(0))=\text{det}\,(S_{q}(\partial B(d_{\partial D}(p)))-S_q(\partial D)).
\tag {3.6}
$$

Using the finiteness of $\Lambda_{\partial D}(p)$ and the partition of the unity, we see that
for the study of the asymptotic expansion of $I_1(\tau)$ and $I_2(\tau)$ it suffices to compute the asymptotic expansion
of the following two integrals, respectively if necessary choosing a smaller $r_q$:
$$\displaystyle
\tilde{I}_1(\tau)
=e^{-2\tau\text{dist}\,(D,B)}\int_{\vert \sigma\vert<r_q}g_0(\sigma)e^{2\tau\phi(\sigma)}d\sigma
$$
and
$$\displaystyle
\tilde{I}_2(\tau)
=e^{-2\tau\text{dist}\,(D,B)}\int_{\vert \sigma\vert<r_q}g_1(p)e^{2\tau\phi(\sigma)}d\sigma,
$$
where
$$\displaystyle
g_0(\sigma)=\Psi_q(\sigma)^{-3}
(\nabla h(\sigma)\cdot\sigma+d_{\partial D}(p)-h(\sigma))
\tag {3.7}
$$
and
$$\displaystyle
g_1(\sigma)
=\frac{g_0(\sigma)}{\Psi_q(\sigma)}-\frac{\beta(x)}{\Psi_q(\sigma)^2}\sqrt{1+\vert\nabla h(\sigma)\vert^2}.
$$
Since we have (3.6), by the Laplace method (\cite{BH}) we obtain
$$\displaystyle
\sqrt{\text{det}\,(\phi_{\sigma_i\sigma_j}(0))}e^{2\tau\text{dist}\,(D,B)}\tilde{I}_1(\tau)
=\frac{\pi}{\tau
d_{\partial D}(p)^2}+O(\tau^{-2}),
\tag {3.8}
$$
$$\displaystyle
\sqrt{\text{det}\,(\phi_{\sigma_i\sigma_j}(0))}e^{2\tau\text{dist}\,(D,B)}\tilde{I}_2(\tau)
=\frac{\pi}{\tau}
\left(\frac{1}{d_{\partial D}(p)^3}-\frac{\beta(q)}{d_{\partial D}(p)^2}\right)
+O(\tau^{-2}).
\tag {3.9}
$$

Thus Theorem 1.1 directly follows from (3.5), (3.8) and (3.9).

For the proof of Theorem 1.2 we have to compute the second term of the asymptotic expansion of $\tilde{I}_1(\tau)$ as
$\tau\longrightarrow\infty$ since the term will make a contribution to the second order term on the right-hand side
of (3.5).
This can be done, in principle, see (8.3.50) on page 338 in \cite{BH}.

However, as pointed out there (see (8.3.53) therein), the coefficient is quite complicated except for the leading term:
$$\begin{array}{c}
\displaystyle
\sqrt{\text{det}\,(\phi_{\sigma_i\sigma_j}(0))}e^{2\tau\text{dist}\,(D,B)}\tilde{I_1}(\tau)
-\frac{\pi}{\tau d_{\partial D}(p)^2}\\
\\
\displaystyle
=
\frac{\pi}{4\tau^2}\triangle_{\xi}G_0\vert_{\xi=0}\sqrt{\text{det}\,(\phi_{\sigma_i\sigma_j}(0))}+O(\tau^{-3}),
\end{array}
\tag {3.10}
$$
where $G_0(\xi)$ is given by (8.3.26) on page 334 in \cite{BH} and $\triangle G_0\vert_{\xi=0}$ has the form in our notation
$$\begin{array}{c}
\displaystyle
\triangle_{\xi}G_0\vert_{\xi=0}
\sqrt{\text{det}\,(\phi_{\sigma_i\sigma_j}(0))}\\
\\
\displaystyle
=\phi_{\sigma_s\sigma_r\sigma_q}(0)B_{sq}B_{rp}(g_0)_{\sigma_p}(0)
-\text{Trace}\,(CB)\\
\\
\displaystyle
-g_0(0)
\left\{\phi_{\sigma_p\sigma_q\sigma_r}(0)\phi_{\sigma_s\sigma_t\sigma_u}(0)
\left(\frac{1}{4}B_{ps}B_{qr}B_{tu}
+\frac{1}{6}B_{ps}B_{qt}B_{ru}\right)
-\frac{1}{4}
\phi_{\sigma_p\sigma_q\sigma_r\sigma_s}(0)B_{pr}B_{qs}\right\},
\end{array}
\tag {3.11}
$$
where
$$\begin{array}{c}
\displaystyle
(B_{pq})=(\phi_{\sigma_p\sigma_q}(0))^{-1},
\\
\\
\displaystyle
C=((g_0)_{\sigma_p\sigma_q}(0)).
\end{array}
$$
Note that we have used the summation convention where repeated
indices are to be summed from $1$ to $2$.

Thus, from (3.9) and (3.10) we obtain
$$\begin{array}{c}
\displaystyle
\sqrt{\text{det}\,(\phi_{\sigma_i\sigma_j}(0))}
e^{2\tau\text{dist}\,(D,B)}\left(\tilde{I}_1(\tau)+\frac{1}{\tau}\tilde{I}_2(\tau)\right)\\
\\
\displaystyle
=\frac{\pi}{d_{\partial D}(p)^2}\frac{1}{\tau}
+\pi
\left(
\frac{1}{4}\triangle_{\xi}G_0\vert_{\xi=0}\sqrt{\text{det}\,(\phi_{\sigma_i\sigma_j}(0))}
-\frac{\beta (q)}{d_{\partial D}(p)^2}+\frac{1}{d_{\partial D}(p)^3}\right)
\frac{1}{\tau^2}\\
\\
\displaystyle
+O\left(\frac{1}{\tau^3}\right)\\
\\
\displaystyle
=\frac{A}{\tau}+\frac{B}{\tau^2}++O\left(\frac{1}{\tau^3}\right),
\end{array}
$$
where
$$\displaystyle
A=\frac{\pi}{d_{\partial D}(p)^2}
$$
and
$$\displaystyle
B=\pi
\left(
\frac{1}{4}\triangle_{\xi}G_0\vert_{\xi=0}\sqrt{\text{det}\,(\phi_{\sigma_i\sigma_j}(0))}
-\frac{\beta (q)}{d_{\partial D}(p)^2}+\frac{1}{d_{\partial D}(p)^3}\right).
$$
Thus, we have
$$\begin{array}{c}
\displaystyle
\sqrt{\text{det}\,(\phi_{\sigma_i\sigma_j}(0))}
e^{2\tau\text{dist}\,(D,B)}\left(\eta-\frac{1}{\tau}\right)^2\left(\tilde{I}_1(\tau)+\frac{1}{\tau}\tilde{I}_2(\tau)\right)\\
\\
\displaystyle
=\frac{\eta^2 A}{\tau}+\frac{\eta}{\tau^2}\left(\eta B-2A\right)+O(\frac{1}{\tau^3}).
\end{array}
\tag {3.12}
$$

Here we prepare the following formulae whose proofs are given in Appendix.
\proclaim{\noindent Lemma 3.1.}  We have
$$\displaystyle
(g_0)_{\sigma_p}(0)=0,
\tag {3.13}
$$
$$\displaystyle
\phi_{\sigma_p\sigma_q\sigma_r}(0)
=h_{\sigma_p\sigma_q\sigma_r}(0),
\tag {3.14}
$$
$$\begin{array}{c}
\displaystyle
\phi_{\sigma_p\sigma_q\sigma_r\sigma_s}(0)B_{pr}B_{qs}
=h_{\sigma_p\sigma_q\sigma_r\sigma_s}(0)B_{pr}B_{qs}
+\frac{14-16dH_{\partial D}(q)}{\displaystyle
d^3
\text{det}\,\phi_{\sigma_i\sigma_j}(0)
}
\end{array}
\tag {3.15}
$$
and
$$\begin{array}{c}
\displaystyle
\text{Trace}\,(CB)
=\frac{8}{d^3}
-\frac{2(1-dH_{\partial D}(q))}{\displaystyle
d^5
\text{det}\,\phi_{\sigma_i\sigma_j}(0)
}.
\end{array}
\tag {3.16}
$$

\endproclaim

From this lemma we obtain
$$\begin{array}{c}
\displaystyle
\triangle_{\xi}G_0\vert_{\xi=0}
\sqrt{\text{det}\,(\phi_{\sigma_i\sigma_j}(0))}
\\
\\
\displaystyle
=-\frac{8}{d^3}+\frac{11-12dH_{\partial D}(q)}
{\displaystyle
2d^5\,\text{det}\,(\phi_{\sigma_i\sigma_j}(0))}\\
\\
\displaystyle
-d^{-2}
h_{\sigma_p\sigma_q\sigma_r}(0)h_{\sigma_s\sigma_t\sigma_u}(0)
\left(\frac{1}{4}B_{ps}B_{qr}B_{tu}+\frac{1}{6}B_{ps}B_{qt}B_{ru}\right)\\
\\
\displaystyle
+\frac{1}{4d^2}
h_{\sigma_p\sigma_q\sigma_r\sigma_s}(0)B_{pr}B_{qs},
\end{array}
\tag {3.17}
$$
where $d=d_{\partial D}(p)$. Thus if one knows the values of all
the second, third and fourth-order derivatives of $h$ at
$\sigma=0$, then the right-hand side on (3.17) is known. Finally
summing up (3.12) over all $q\in\Lambda_{\partial D}(p)$,
from (3.5) we see that (1.6) in Theorem 1.2 is true.
The proof of Lemma
3.1 is given in Appendix.

\section{A data reduction and its implication}

\subsection{Proof of Proposition 1.1.}

Let $U$ be a non-empty open set of $\Bbb R^3$.
The mean value theorem \cite{CH} states that if $\Psi\in C^2(U)$ is a solution of the modified Helmholtz equation
$$\displaystyle
(\triangle-\tau^2)\Psi=0\,\,\text{in}\,U
$$
and $\overline{B_r(x)}\subset U$, then
$$\displaystyle
\Psi(x)
=\frac{\tau}{4\pi r
\sinh (\tau r)}
\int_{\partial B_r(x)}\Psi(y)dS(y).
\tag {4.1}
$$
Substituting this into the identity
$$\displaystyle
\int_{B_r(x)}\Psi(z)dz=\int_0^r\left(\int_{\partial B_s(x)}\Psi(y)dS(y)\right)ds,
$$
we obtain
$$\displaystyle
\int_{B_r(x)}\Psi(z)dz
=\frac{4\pi}{\tau^3}
(\tau r\cosh(\tau r)-\sinh(\tau r))\Psi(x).
\tag {4.2}
$$

\proclaim{\noindent Lemma 4.1.}
Let $B$ and $f$ be same as those of Problem I.
Let $B_R(p)$ denote the open ball centred at $p$ with radius $R$
and satisfy $B_R(p)\subset\Bbb R^3\setminus\overline D$ and $B\subset B_R(p)$.
We have
$$\begin{array}{c}
\displaystyle
\int_{\partial B_R(p)}(w_f-v_f)dS
\\
\\
\displaystyle
=\frac{\tau^2 R\sinh(\tau R)}{\tau\eta\cosh\,(\tau\eta)-\sinh\,(\tau\eta)}
\int_B(w_f-v_f)dx+O(\tau e^{-\tau(T-(R-\eta))}).
\end{array}
\tag {4.3}
$$

\endproclaim

{\it\noindent Proof.}
We have
$$\begin{array}{c}
\displaystyle
\int_B(w_f-v_f)dx
=\int_B\epsilon_f^0dx+\int_BZdx.
\end{array}
$$
Since we have (2.5), it holds that
$$\displaystyle
\int_BZdx=O(e^{-\tau T})
$$
and thus
$$\begin{array}{c}
\displaystyle
\int_B(w_f-v_f)dx
=\int_B\epsilon_f^0dx+O(e^{-\tau T}).
\end{array}
\tag {4.4}
$$

Applying (4.1) and (4.2) to $\epsilon_f^0$ and putting
$\varphi(\xi)=\xi\cosh\xi-\sinh\xi$, we have
$$\begin{array}{c}
\displaystyle
\int_B\epsilon_f^0 dx
=\frac{4\pi}{\tau^3}\varphi(\tau\eta)\epsilon_f^0(p)
\\
\\
\displaystyle
=\frac{4\pi}{\tau^3}\varphi(\tau\eta)\times
\frac{\tau}{4\pi R\sinh\,(\tau R)}\int_{\partial B_R(p)}\epsilon_f^0 dS
\\
\\
\displaystyle
=\frac{\varphi(\tau\eta)}{\tau^2 R\sinh(\tau R)}\int_{\partial B_R(p)}(w_f-v_f-Z)dS\\
\\
\displaystyle
=\frac{\varphi(\tau\eta)}{\tau^2 R\sinh(\tau R)}
\left(\int_{\partial B_R(p)}(w_f-v_f)dS
-\int_{\partial B_R(p)}ZdS\right).
\end{array}
$$

Just simply applying the trace theorem onto $\partial B_R(p)$, from (4.2) we have
$$\displaystyle
\int_{\partial B_R(p)}ZdS=O(\tau e^{-\tau T}).
$$
Thus, we obtain
$$\begin{array}{c}
\displaystyle
\int_B\epsilon_f^0dx
=\frac{\varphi(\tau\eta)}{\tau^2 R\sinh(\tau R)}
\left(\int_{\partial B_R(p)}(w_f-v_f)dS+O(\tau e^{-\tau T})\right).
\end{array}
$$
From this and (4.4) we obtain
$$\begin{array}{c}
\displaystyle
\int_B(w_f-v_f)dx
=\frac{\varphi(\tau\eta)}{\tau^2 R\sinh(\tau R)}\left(\int_{\partial B_R(p)}(w_f-v_f)dS
+O(\tau e^{-\tau T})\right)
+O(e^{-\tau T}).
\end{array}
\tag {4.5}
$$
Since $R\ge\eta$ and
$$\displaystyle
O(e^{-\tau T}\frac{\tau^2 R\sinh(\tau R)}{\varphi(\tau\eta)})
=O(e^{-\tau T}\frac{\tau^2 Re^{\tau R}}{\tau\eta e^{\tau\eta}})
=O(\tau e^{-\tau T}e^{\tau(R-\eta)}),
$$
from (4.5) one gets (4.3).

\noindent
$\Box$

Since we have
$$\begin{array}{c}
\displaystyle
\frac{\tau^2 R\sinh(\tau R)}{\tau\eta\cosh\,(\tau\eta)-\sinh\,(\tau\eta)}\\
\\
\displaystyle
=\frac{\tau^2 R(e^{\tau R}-e^{-\tau R})}
{\tau\eta(e^{\tau\eta}+e^{-\tau\eta})-(e^{\tau\eta}-e^{-\tau\eta})}\\
\\
\displaystyle
=\frac{\tau^2 Re^{\tau R}(1-e^{-2\tau R})}
{\tau\eta e^{\tau\eta}(1+e^{-2\tau\eta})-e^{\tau\eta}(1-e^{-2\tau\eta})}\\
\\
\displaystyle
=\frac{\tau^2 Re^{\tau R}(1-e^{-2\tau R})}
{\tau\eta e^{\tau\eta}(1+O(\tau^{-1}e^{-\tau\eta}))}\\
\\
\displaystyle
=\frac{R}{\eta}\tau e^{\tau(R-\eta)}
\frac{1-e^{-2\tau R}}{1+O(\tau^{-1}e^{-\tau\eta})}\\
\\
\displaystyle
=\frac{R}{\eta}\tau e^{\tau(R-\eta)}
\left(\frac{1}{1+O(\tau^{-1}e^{-\tau\eta})}+O(e^{-2\tau R})\right)\\
\\
\displaystyle
=\frac{R}{\eta}\tau e^{\tau(R-\eta)}
(1+O(\tau^{-1}e^{-\tau\eta})+O(e^{-2\tau R}))\\
\\
\displaystyle
=\frac{R}{\eta}\tau e^{\tau(R-\eta)}(1+O(\tau^{-1}e^{-\tau\eta})),
\end{array}
$$
(4.3) yields (1.7).  This completes the proof of Proposition 1.1.

\subsection{Some transplanted results}

Let $B$ and $f$ be same as those of Problem I.
Let $B_R(p)$ denote the open ball centred at $p$ with radius $R$
and satisfy $B_R(p)\subset\Bbb R^3\setminus\overline D$ and $B\subset B_R(p)$.

Let $T>2\text{dist}\,(D,B)$.
From (1.7) we have.
$$\begin{array}{c}
\displaystyle
\frac{\tau^3\eta}{R}e^{\tau(2\text{dist}\,(D,B)-(R-\eta))}
\int_{\partial B_R(p)}(w_f-v_f)dS
\\
\\
\displaystyle
=\tau^4 e^{2\tau\text{dist}\,(D,B)}
(1+O(\tau^{-1}e^{-\tau\eta}))\int_B(w_f-v_f)dx
+O(\tau^4e^{-\tau(T-2\text{dist}\,(D,B))}).
\end{array}
\tag {4.6}
$$
This enables us to transplant several results with data $u_f$ on
$B\times\,]0,\,T[$ provided $T>2\text{dist}\,(D,B)$ not
$T>2\text{dist}\,(D,B)-(R-\eta)$.  We think this is reasonable
since the quantity $2\text{dist}\,(D,B)-(R-\eta)$ coincides with
the time of flight of wave with propagation speed $1$ that starts
at points on $\partial B$ and hits points on $\partial D$ and
returns to $\partial B_R(p)$.  However, the surface is {\it thin},
so to make the received signal on $\partial B_R(p)$ strong enough
some redundant observation time may be needed.

Using (4.6), from Theorem 1.2 in \cite{IE3} one can deduce the following result.

\proclaim{\noindent Corollary 4.1.}
Assume that $\partial D$ is $C^2$.
Let $T$ satisfy $T>2\text{dist}\,(D,B)$.

We have the following:

if $\gamma(x)\le 1-C$ a.e. $x\in\partial D$, then there exists $\tau_0>0$ such that, for all $\tau\ge\tau_0$
$$
\displaystyle
\int_{\partial B_R(p)}(w_f-v_f)dS>0;
$$

if $\gamma(x)\ge 1+C$ a.e. $x\in\partial D$ for a positive constant $C$, then there exists $\tau_0>0$ such that, for all $\tau\ge\tau_0$
$$
\displaystyle
\int_{\partial B_R(p)}(w_f-v_f)dS<0.
$$

In both cases, the formula
$$\displaystyle
\lim_{\tau\longrightarrow\infty}\frac{1}{\tau}\log\left\vert\int_{\partial B_R(p)}(w_f-v_f)dS\right\vert=
-(2\text{dist}\,(D,B)-(R-\eta)),
$$
is valid.
\endproclaim

The following result is a transplanted version of Theorems 1.1 via (4.6).

\proclaim{\noindent Corollary 4.2.}  Let $\gamma\equiv 0$.
Assume that $\partial D$ is $C^3$ and $\beta\in C^2(\partial D)$; $\Lambda_{\partial D}(p)$ is finite and satisfies
$$\displaystyle
\text{det}\,(S_q(\partial B_{d_{\partial D}(p)}(p))-S_q(\partial D))>0,\,\,\forall q\in\Lambda_{\partial D}(p).
$$
If $T$ satisfies $T>2\text{dist}\,(D,B)$,
then we have
$$\begin{array}{c}
\displaystyle
\lim_{\tau\longrightarrow\infty}
\tau^3\frac{\eta}{R}e^{\tau(2\text{dist}\,(D,B)-(R-\eta))}
\int_{\partial B_R(p)}(w_f-v_f)dS
=\frac{\pi}{2}
\left(\frac{\eta}{d_{\partial D}(p)}\right)^2A_{\partial D}(p).
\end{array}
$$

\endproclaim

From (4.6) we have
$$\begin{array}{c}
\displaystyle
\tau^4
\left\{
\frac{\eta}{R} e^{\tau(2\text{dist}\,(D,B)-(R-\eta))}\int_{\partial B_R(p)}(w_f-v_f)dS
-\frac{1}{\tau^3}\frac{\pi}{2}\left(\frac{\eta}{d_{\partial D}(p)}\right)^2A_{\partial D}(p)\right\}\\
\\
\displaystyle
=(1+O(\tau^{-1}e^{-\tau\eta}))\tau^5
\left\{
e^{2\tau\text{dist}\,(D,B)}\int_{\partial B}(w_f-v_f)dx
-\frac{1}{\tau^4}\frac{\pi}{2}\left(\frac{\eta}{d_{\partial D}(p)}\right)^2A_{\partial D}(p)\right\}
\\
\\
\displaystyle
+O(e^{-\tau\eta})+O(\tau^5 e^{-\tau(T-2\text{dist}\,(D,B))}).
\end{array}
$$

This enables us to deduce the following result from Theorem 1.2.

\proclaim{\noindent Corollary 4.3.}  Let $\gamma\equiv 0$.
Assume that $\partial D$ is $C^5$ and $\beta\in C^2(\partial D)$;
$\Lambda_{\partial D}(p)$ is finite and satisfies
$$\displaystyle
\text{det}\,(S_q(\partial B_{d_{\partial D}(p)}(p))-S_q(\partial D))>0,\,\,\forall q\in\Lambda_{\partial D}(p).
$$
For each $q\in \Lambda_{\partial D}(p)$ let $\mbox{\boldmath
$e$}_j$, $j=1,2$ be orthogonal basis of the tangent space at $q$
of $\partial D$ with $\mbox{\boldmath $e$}_1\times\mbox{\boldmath
$e$}_2=\nu_q$. Choose an open ball $U$ centred at $q$ with radius
$r_q$ in such a way that: there exist a $h\in C^5_0(\Bbb R^2)$
with $h(0,0)=0$ and $\nabla h(0,0)=0$ such that $U\cap\partial
D=\{q+\sigma_1\mbox{\boldmath $e$}_1+\sigma_2\mbox{\boldmath
$e$}_2+h(\sigma_1,\sigma_2)\nu_q\,\vert\,
\sigma_1^2+\sigma_2^2+h(\sigma_1,\sigma_2)^2<r_q^2\}$.

If $T$ satisfies $T>2\text{dist}\,(D,B)$,
then we have
$$\begin{array}{c}
\displaystyle
\lim_{\tau\longrightarrow\infty}
\tau^4\left\{
e^{\tau(2\text{dist}\,(D,B)-(R-\eta))}\frac{\eta}{R}
\int_{\partial B_R(p)}(w_f-v_f)dS
-\frac{1}{\tau^3}\frac{\pi}{2}
\left(\frac{\eta}{d_{\partial D}(p)}\right)^2A_{\partial D}(p)\right\}\\
\\
\displaystyle
=-\frac{\pi\eta}{d_{\partial D}(p)^2}A_{\partial D}(p)
+\frac{\pi}{2}\eta^2 B_{\partial D}(p).
\end{array}
$$

\endproclaim

Some trivial modifications of Corollaries 1.1 and 1.2 are also
valid, however, we omit their description.

\section{Concluding remarks}

In this paper we considered only the case when $\gamma\equiv 0$.
Thus, the next problem is to consider general $\gamma(\ge 0)$.
The author thinks that Proposition 2.1 should be changed because of the {\it dissipation}
of the energy of the wave on the obstacle.

In \cite{IE5} we have considered the case when the so-called {\it bistatic data} are available.
It means that one observes the wave at a different place, for example,
even possible being far from the center of the support of an initial data.
We found that a {\it spheroid} plays a role similar to a sphere.
However, the technique used in \cite{IE5} based on the maximum principle and
depends heavily on the homogeneous Dirichlet boundary condition on the obstacle
which corresponds to the case when $\gamma\equiv 0$ and $\beta=\infty$.
Thus, one can not be readily apply the method to the present case.

For further open problems see section 6 of \cite{IE5} and also section 4 of \cite{IE4}.

$$\quad$$

\centerline{{\bf Acknowledgement}}

This research was partially supported by Grant-in-Aid for
Scientific Research (C)(No. 25400155) of Japan  Society for the
Promotion of Science.

\section{Appendix}

\subsection{A direct derivation of (3.4)}

Applying the mean value theorem (4.2) to $v_f(x)$ for $x\in\Bbb R^3\setminus\overline B$ with $B=B_{\eta}(p)$, we obtain
$$
\displaystyle
v_f(x)
=\frac{\varphi(\tau\eta)}{\tau^3}\frac{e^{-\tau\vert x-p\vert}}{\vert x-p\vert}
$$
and thus
$$\displaystyle
\nabla v_f(x)
=\frac{\varphi(\tau\eta)}{\tau^2}
\left(1+\frac{1}{\tau\vert x-p\vert}\right)\frac{(p-x)}{\vert x-p\vert}\frac{e^{-\tau\vert x-p\vert}}{\vert x-p\vert},
$$
where
$$\displaystyle
\varphi(\xi)=\xi\cosh\xi-\sinh\xi.
$$
Substituting these into $J(\tau)$, we obtain
$$\begin{array}{c}
\displaystyle
J(\tau)=\frac{\varphi(\tau\eta)^2}{\tau^5}
\left(I_1(\tau)+\frac{1}{\tau}I_2(\tau)\right),
\end{array}
$$
where
$$\displaystyle
\tilde{I_1}(\tau)
=\int_{\partial D}\frac{1}{\vert x-p\vert^2}\frac{(p-x)}{\vert x-p\vert}\cdot\nu_x\,
e^{-2\tau\vert x-p\vert}
dS_x
$$
and
$$\displaystyle
\tilde{I_2}(\tau)
=\int_{\partial D}
\left(\frac{1}{\vert x-p\vert^3}\frac{(p-x)}{\vert x-p\vert}\cdot\nu_x
-\frac{\beta (x)}{\vert x-p\vert^2}\right)e^{-2\tau\vert x-p\vert}
dS_x.
$$
Since we have
$$\varphi(\xi)e^{-\xi}
=\frac{\xi-1}{2}+\frac{\xi+1}{2}e^{-2\xi},
$$
one gets
$$\displaystyle
\frac{\varphi(\tau\eta)^2e^{-2\tau\eta}}{\tau^5}
=\frac{1}{4\tau^3}
\left(\eta-\frac{1}{\tau}\right)^2
+O(\tau^{-3}e^{-2\tau\eta}).
$$
Note that also $d_B(x)=\vert x-p\vert-\eta$ for $x\in\partial D$
and thus, $I_j(\tau)=e^{2\tau\eta}\tilde{I}_j(\tau)$, $j=1,2$.

From these one obtains (3.4).

\subsection{Proof of Lemma 3.1}

In this section we set $d=d_{\partial D}(p)$.

\subsubsection{Proof of (3.13)}

From (3.7)
we have
$$\begin{array}{c}
\displaystyle
\frac{\partial}{\partial\sigma_p}g_0
=-3\Psi^{-4}\frac{\partial\Psi}{\partial\sigma_p}(h_{\sigma_q}\sigma_q+d-h)
+\Psi^{-3}
\left(h_{\sigma_p\sigma_q}\sigma_q+h_{\sigma_p}-h_{\sigma_p}\right)\\
\\
\displaystyle
=-3\Psi^{-1}\frac{\partial\Psi}{\partial\sigma_p}g_0
+\Psi^{-3}h_{\sigma_p\sigma_q}\sigma_q.
\end{array}
\tag {A.1}
$$
Thus we obtain (3.13).

\subsubsection{Proof of (3.14)}

We have
$$\displaystyle
\frac{\partial\Psi_q}{\partial\sigma_p}(\sigma)
=\Psi_q(\sigma)^{-1}
\{\sigma_p+h_{\sigma_p}(\sigma)(h(\sigma)-d)\},
\tag {A.2}
$$
This gives
$$\displaystyle
\phi_{\sigma_p}(\sigma)
=-\Psi_q(\sigma)^{-1}
\{\sigma_p+h_{\sigma_p}(\sigma)(h(\sigma)-d)\}
\tag {A.3}
$$
and
$$\displaystyle
\phi_{\sigma_p}(0)=\frac{\partial \Psi_{q}}{\partial\sigma_p}(0)=0.
\tag {A.4}
$$
It follows from (A.3) that
$$\begin{array}{c}
\displaystyle
\phi_{\sigma_p\sigma_q}(\sigma)
=(\Psi_q(\sigma))^{-2}\frac{\partial\Psi_q}{\partial\sigma_q}(\sigma)
\{\sigma_q+h_{\sigma_q}(\sigma)(h(\sigma)-d)\}\\
\\
\displaystyle
-\Psi_q(\sigma)^{-1}
\left\{
\delta_{pq}+h_{\sigma_p\sigma_q}(\sigma)(h(\sigma)-d)
+h_{\sigma_p}(\sigma)h_{\sigma_q}(\sigma)\right\}.
\end{array}
$$
Applying (A.2) to the first term in this right-hand side, we obtain
$$\begin{array}{c}
\displaystyle
\Psi_q(\sigma)\phi_{\sigma_p\sigma_q}(\sigma)
=
\phi_{\sigma_p}(\sigma)\phi_{\sigma_q}(\sigma)
-\left\{
\delta_{pq}+h_{\sigma_p\sigma_q}(\sigma)(h(\sigma)-d)
+h_{\sigma_p}(\sigma)h_{\sigma_q}(\sigma)\right\}.
\end{array}
\tag {A.5}
$$
Differentiating both sides on (A.5) with respect to $\sigma_r$ and using (A.2) and (A.3), we obtain
$$\begin{array}{c}
\displaystyle
-\phi_{\sigma_r}(\sigma)\phi_{\sigma_p\sigma_q}(\sigma)
+\Psi_q(\sigma)\phi_{\sigma_p\sigma_q\sigma_r}(\sigma)
=
\phi_{\sigma_p\sigma_r}(\sigma)\phi_{\sigma_q}(\sigma)
+\phi_{\sigma_p}(\sigma)\phi_{\sigma_q\sigma_r}(\sigma)\\
\\
\displaystyle
-\left\{
h_{\sigma_p\sigma_q\sigma_r}(\sigma)(h(\sigma)-d)
+h_{\sigma_p\sigma_q}(\sigma)h_{\sigma_r}(\sigma)
+
h_{\sigma_p\sigma_r}(\sigma)h_{\sigma_q}(\sigma)
+h_{\sigma_p}(\sigma)h_{\sigma_q\sigma_r}(\sigma)
\right\}.
\end{array}
\tag {A.6}
$$
Noting $h(0)=h_{\sigma_p}(0)=0$, from this together with (A.4) we obtain (3.14).

\subsubsection{Proof of (3.15)}

From (A.4) we have
$$\begin{array}{c}
\displaystyle
\frac{\partial}{\partial\sigma_s}
\left(-\phi_{\sigma_r}(\sigma)\phi_{\sigma_p\sigma_q}(\sigma)
+\Psi_q(\sigma)\phi_{\sigma_p\sigma_q\sigma_r}(\sigma)\right)\vert_{\sigma=0}\\
\\
\displaystyle
=-\phi_{\sigma_r\sigma_s}(0)\phi_{\sigma_p\sigma_q}(0)
+d\,\phi_{\sigma_p\sigma_q\sigma_r\sigma_s}(0)
\end{array}
\tag {A.7}
$$
and also
$$\begin{array}{c}
\displaystyle
\frac{\partial}{\partial\sigma_s}\left(\phi_{\sigma_p\sigma_r}(\sigma)\phi_{\sigma_q}(\sigma)
+\phi_{\sigma_p}(\sigma)\phi_{\sigma_q\sigma_r}(\sigma)\right)\vert_{\sigma=0}\\
\\
\displaystyle
=\phi_{\sigma_p\sigma_r}(0)\phi_{\sigma_q\sigma_s}(0)
+\phi_{\sigma_p\sigma_s}(0)\phi_{\sigma_q\sigma_r}(0).
\end{array}
\tag {A.8}
$$
Moreover we have
$$\begin{array}{c}
\displaystyle
\frac{\partial}{\partial\sigma_s}\left\{
h_{\sigma_p\sigma_q\sigma_r}(\sigma)(h(\sigma)-d)
+h_{\sigma_p\sigma_q}(\sigma)h_{\sigma_r}(\sigma)
+
h_{\sigma_p\sigma_r}(\sigma)h_{\sigma_q}(\sigma)
+h_{\sigma_p}(\sigma)h_{\sigma_q\sigma_r}(\sigma)
\right\}\vert_{\sigma=0}\\
\\
\displaystyle
=
-d\,h_{\sigma_p\sigma_q\sigma_r\sigma_s}(0)\\
\\
\displaystyle
+h_{\sigma_p\sigma_q}(0)h_{\sigma_r\sigma_s}(0)
+
h_{\sigma_p\sigma_r}(0)h_{\sigma_q\sigma_s}(0)
+h_{\sigma_p\sigma_s}(0)h_{\sigma_q\sigma_r}(0).
\end{array}
\tag {A.9}
$$
Differentiating both sides on (A.6) with respect to $\sigma_s$ and using (A.7)-(A.9), we obtain
$$\begin{array}{c}
\displaystyle
\phi_{\sigma_p\sigma_q\sigma_r\sigma_s}(0)
=\frac{1}{d}\left(\phi_{\sigma_r\sigma_s}(0)\phi_{\sigma_p\sigma_q}(0)-h_{\sigma_p\sigma_q}(0)h_{\sigma_r\sigma_s}(0)\right)\\
\\
\displaystyle
+\frac{1}{d}
\left(
\phi_{\sigma_p\sigma_r}(0)\phi_{\sigma_q\sigma_s}(0)-h_{\sigma_p\sigma_r}(0)h_{\sigma_q\sigma_s}(0)\right)\\
\\
\displaystyle
+\frac{1}{d}
\left(\phi_{\sigma_p\sigma_s}(0)\phi_{\sigma_q\sigma_r}(0)-h_{\sigma_p\sigma_s}(0)h_{\sigma_q\sigma_r}(0)\right)
\\
\\
\displaystyle
+\,h_{\sigma_p\sigma_q\sigma_r\sigma_s}(0).
\end{array}
$$
A combination of this and (A.5) yields
$$\begin{array}{c}
\displaystyle
\phi_{\sigma_p\sigma_q\sigma_r\sigma_s}(0)
=\frac{1}{d^3}\left(\delta_{pq}\delta_{rs}+\delta_{qr}\delta_{ps}+\delta_{qs}\delta_{pr}\right)
\\
\\
\displaystyle
-\frac{1}{d^2}
\left(
\delta_{pq}h_{\sigma_r\sigma_s}(0)
+\delta_{rs}h_{\sigma_p\sigma_q}(0)
+\delta_{qr}h_{\sigma_p\sigma_s}(0)
+\delta_{ps}h_{\sigma_q\sigma_r}(0)
+\delta_{qs}h_{\sigma_p\sigma_r}(0)
+\delta_{pr}h_{\sigma_q\sigma_s}(0)
\right)
\\
\\
\displaystyle
+h_{\sigma_p\sigma_q\sigma_r\sigma_s}(0).
\end{array}
\tag {A.10}
$$
We have also
$$\begin{array}{c}
\displaystyle
\delta_{pq}\delta_{rs}B_{pr}B_{qs}
=B_{qs}B_{qs},\,
\delta_{qr}\delta_{ps}B_{pr}B_{qs}
=B_{pq}B_{qp},\,
\delta_{qs}\delta_{pr}B_{pr}B_{qs}
=B_{pp}B_{qq}.
\end{array}
$$
Since $B_{pq}=B_{qp}$ and $B_{pp}=\text{Trace}\,B$, from these
we obtain
$$\displaystyle
\left(\delta_{pq}\delta_{rs}+\delta_{qr}\delta_{ps}+\delta_{qs}\delta_{pr}\right)B_{pr}B_{qs}
=2\vert B\vert^2+(\text{Trace}\,B)^2.
\tag {A.11}
$$

Moreover, we have
$$\begin{array}{c}
\displaystyle
\delta_{pq}h_{\sigma_r\sigma_s}(0)B_{pr}B_{qs}
=\text{Trace}\,(B\nabla^2h(0)B);
\\
\\
\displaystyle
\delta_{rs}h_{\sigma_p\sigma_q}(0)B_{pr}B_{qs}
=\text{Trace}\,(\nabla^2h(0)B^2);
\\
\\
\displaystyle
\delta_{qr}h_{\sigma_p\sigma_s}(0)B_{pr}B_{qs}
=\text{Trace}\,(\nabla^2h(0)B^2);
\\
\\
\displaystyle
\delta_{ps}h_{\sigma_q\sigma_r}(0)B_{pr}B_{qs}
=\text{Trace}\,(\nabla^2h(0)B^2);
\\
\\
\displaystyle
\delta_{qs}h_{\sigma_p\sigma_r}(0)B_{pr}B_{qs}
=\text{Trace}\,B\,\text{Trace}\,(\nabla^2h(0)B);
\\
\\
\displaystyle
\delta_{pr}h_{\sigma_q\sigma_s}(0)B_{pr}B_{qs}
=\text{Trace}\,B\,\text{Trace}\,(\nabla^2h(0)B).
\end{array}
$$
Now it follows from these, (A.10) and (A.11) that
$$\begin{array}{c}
\displaystyle
\phi_{\sigma_p\sigma_q\sigma_r\sigma_s}(0)B_{pr}B_{qs}
=\frac{1}{d^3}\{2\vert B\vert^2+(\text{Trace}\,B)^2\}
\\
\\
\displaystyle
-\frac{1}{d^2}
\left(\text{Trace}\,(B\nabla^2h(0)B)
+3\text{Trace}\,(\nabla^2h(0)B^2)
+2\text{Trace}\,B\,\text{Trace}\,(\nabla^2h(0)B)
\right)
\\
\\
\displaystyle
+h_{\sigma_p\sigma_q\sigma_r\sigma_s}(0)B_{pr}B_{qs}.
\end{array}
\tag {A.12}
$$
Note that
$$\begin{array}{c}
\displaystyle
\nabla^2h(0)=B^{-1}+\frac{1}{d}I_2;\\
\\
\displaystyle
B\nabla^2h(0)B=\nabla^2h(0)B^2
=B+\frac{1}{d}B^2;
\\
\\
\displaystyle
\text{Trace}\,(\nabla^2h(0)B)
=2+\frac{1}{d}\text{Trace}\,B;\\
\\
\displaystyle
\text{Trace}\,(B\nabla^2h(0)B)
=\text{Trace}\,(\nabla^2h(0)B^2)
=\text{Trace}\,B+\frac{1}{d}\vert B\vert^2.
\end{array}
\tag {A.13}
$$
Substituting these into the right-hand side on (A.12), we obtain
$$\begin{array}{c}
\displaystyle
d^3(\phi_{\sigma_p\sigma_q\sigma_r\sigma_s}(0)B_{pr}B_{qs}-h_{\sigma_p\sigma_q\sigma_r\sigma_s}(0)B_{pr}B_{qs})
\\
\\
\displaystyle
=\vert B\vert^2-(\text{Trace}\,B)^2-8d\,\text{Trace}\,B.
\end{array}
\tag {A.14}
$$
A direct computation yields
$$\begin{array}{c}
\vert B\vert^2
=
\frac{\displaystyle
\left(\frac{1}{d}-h_{\sigma_1\sigma_1}(0)\right)^2
+\left(\frac{1}{d}-h_{\sigma_2\sigma_2}(0)\right)^2
+2h_{\sigma_1\sigma_2}(0)^2}
{\displaystyle
(\text{det}\,\phi_{\sigma_i\sigma_j}(0))^2},
\\
\\
\displaystyle
\text{Trace}\,B
=-
\frac{
\displaystyle
\left(\frac{1}{d}-h_{\sigma_1\sigma_1}(0)\right)
+
\left(\frac{1}{d}-h_{\sigma_2\sigma_2}(0)\right)
}
{\displaystyle
\text{det}\,\phi_{\sigma_i\sigma_j}(0)}
\end{array}
\tag {A.15}
$$
and thus
$$\begin{array}{c}
\displaystyle
(\text{det}\,\phi_{\sigma_i\sigma_j}(0))^2
\left\{\vert B\vert^2-(\text{Trace}\,B)^2\right\}
\\
\\
\displaystyle
=\left(\frac{1}{d}-h_{\sigma_1\sigma_1}(0)\right)^2
+\left(\frac{1}{d}-h_{\sigma_2\sigma_2}(0)\right)^2
+2h_{\sigma_1\sigma_2}(0)^2
-\left\{\left(\frac{1}{d}-h_{\sigma_1\sigma_1}(0)\right)
+
\left(\frac{1}{d}-h_{\sigma_2\sigma_2}(0)\right)\right\}^2\\
\\
\displaystyle
=2\left\{h_{\sigma_1\sigma_2}(0)^2
-\left(\frac{1}{d}-h_{\sigma_1\sigma_1}(0)\right)
\left(\frac{1}{d}-h_{\sigma_2\sigma_2}(0)\right)\right\}\\
\\
\displaystyle
=-2
\left\vert
\begin{array}{lr}
\displaystyle
\frac{1}{d}-h_{\sigma_1\sigma_1}(0) &
\displaystyle
-h_{\sigma_1\sigma_2}(0)\\
\\
\displaystyle
-h_{\sigma_1\sigma_2}(0) &
\displaystyle
\frac{1}{d}-h_{\sigma_2\sigma_2}(0)
\end{array}
\right\vert
=-2\text{det}\,\phi_{\sigma_i\sigma_j}(0).
\end{array}
$$
This yields
$$\displaystyle
\text{det}\,\phi_{\sigma_i\sigma_j}(0)
\left\{\vert B\vert^2-(\text{Trace}\,B)^2\right\}
=-2.
$$
Then we have
$$\begin{array}{c}
\displaystyle
\text{det}\,\phi_{\sigma_i\sigma_j}(0)
\left\{\vert B\vert^2-(\text{Trace}\,B)^2
-8d\text{Trace}\,B\right\}
=14-16dH_{\partial D}(q).
\end{array}
$$
Now from (A.14) we obtain (3.15).

\subsubsection{Proof of (3.16)}

From (A.1) we have
$$\begin{array}{c}
\displaystyle
\frac{\partial^2}{\partial\sigma_q\partial\sigma_p}g_0
=\frac{\partial}{\partial\sigma_q}
\left(
-3\Psi^{-1}\frac{\partial\Psi}{\partial\sigma_p}g_0
+\Psi^{-3}h_{\sigma_p\sigma_r}\sigma_r
\right)\\
\\
\displaystyle
=
3\Psi^{-2}\frac{\partial\Psi}{\partial\sigma_q}\frac{\partial\Psi}{\partial\sigma_p}g_0
-3\Psi^{-1}\frac{\partial^2\Psi}{\partial\sigma_q\partial\sigma_p}g_0
-3\Psi^{-1}\frac{\partial\Psi}{\partial\sigma_p}\frac{\partial g_0}{\partial\sigma_q}\\
\\
\displaystyle
-3\Psi^{-4}\frac{\partial\Psi}{\partial\sigma_q}h_{\sigma_p\sigma_r}\sigma_r
+\Psi^{-3}h_{\sigma_p\sigma_q\sigma_r}\sigma_r
+\Psi^{-3}h_{\sigma_p\sigma_q}.
\end{array}
$$
This together with (A.4) and (3.13) yields
$$\begin{array}{c}
\displaystyle
(g_0)_{\sigma_p\sigma_q}(0)
=
-3d^{-3}\frac{\partial^2\Psi}{\partial\sigma_q\partial\sigma_p}(0)
+d^{-3}h_{\sigma_p\sigma_q}(0)\\
\\
\displaystyle
=d^{-3}
\left(h_{\sigma_p\sigma_q}(0)-3\frac{\partial^2\Psi}{\partial\sigma_q\partial\sigma_p}(0)
\right).
\end{array}
\tag {A.16}
$$
Since we have
$$\displaystyle
\frac{\partial^2\Psi}{\partial\sigma_p\partial\sigma_q}(0)
=\frac{1}{d}\delta_{pq}-h_{\sigma_p\sigma_q}(0),
$$
from (A.16), we obtain
$$
\displaystyle
(g_0)_{\sigma_p\sigma_q}(0)
=d^{-3}
\left(4h_{\sigma_p\sigma_q}(0)-\frac{3}{d}\delta_{pq}
\right).
$$
This together with (A.13) yields
$$\begin{array}{c}
\displaystyle
\text{Trace}\,(CB)
=(g_0)_{\sigma_p\sigma_q}B_{qp}
\\
\\
\displaystyle
=d^{-3}
\left(4h_{\sigma_p\sigma_q}(0)B_{qp}
-\frac{3}{d}\delta_{pq}B_{qp}\right)\\
\\
\displaystyle
=d^{-3}
\left(
4\text{Trace}\,(\nabla^2h(0)B)
-\frac{3}{d}\text{Trace}\,B\right)\\
\\
\displaystyle
=d^{-3}
\left\{4\left(2+\frac{1}{d}\text{Trace}\,B\right)
-\frac{3}{d}\text{Trace}\,B
\right\}\\
\\
\displaystyle
=d^{-3}
\left(8
+\frac{1}{d}\text{Trace}\,B
\right).
\end{array}
$$
Now from this and the second formula on (A.15) we obtain (3.16).

\vskip1cm
\noindent
e-mail address

ikehata@amath.hiroshima-u.ac.jp


\begin{thebibliography}{99}





\bibitem{BH}  Bleistein, N. and Handelsman, R. A., Asymptotic expansions of integrals, Dover Publications, New York, 1986.





\bibitem{CH} Courant, R. and Hilbert, D.,
              Methoden der Mathematischen Physik, vol. 2., Berlin, Springer, 1937.


\bibitem{DL}  Dautray, R. and Lions, J-L., Mathematical analysis and numerical methods for
              sciences and technology,
              Evolution problems I, Vol. {\bf 5}, Springer-Verlag, Berlin, 1992.




\bibitem{GT}  Gilbarg, D. and Trudinger, N. S., Elliptic partial differential equations
              of second order, second.ed., Springer-Verlag, Berlin, Heidelberg,
              New York,Tokyo, 1983.








\bibitem{I1} Ikehata, M.,
             \newblock Enclosing a polygonal cavity in a two-dimensional bounded domain from Cauchy data,
             \newblock Inverse Problems, {\bf 15}(1999), 1231-1241.




\bibitem{IF} Ikehata, M., Reconstruction of the support function for inclusion from boundary measurements,
             \newblock J. Inv. Ill-Posed Problems, {\bf 8}(2000), 367-378.





\bibitem{ICAL} Ikehata, M., On reconstruction in the inverse conductivity problem
                with one measurement, Inverse Problems, {\bf 16}(2000), 785-793.



\bibitem{IR}  Ikehata, M., A regularized extraction formula in the enclosure method,
              Inverse Problems, {\bf 18}(2002), 435-440.


\bibitem{ISC} Ikehata, M., Inverse scattering problems and the enclosure method,
              Inverse Problems, {\bf 20}(2004), 533-551.



\bibitem{IH}  Ikehata, M., The Herglotz wave function, the Vekua transform and the enclosure method,
              Hiroshima Math. J., {\bf 35}(2005), 485-506.




\bibitem{I4} Ikehata, M.,
             \newblock Extracting discontinuity in a heat conductive body. One-space dimensional case,
             \newblock Applicable Analysis, {\bf 86}(2007), no. 8, 963-1005.





\bibitem{IE0}  Ikehata, M.,
              The enclosure method for inverse obstacle scattering problems with dynamical data  over a finite
              time interval, Inverse Problems, {\bf 26}(2010) 055010(20pp).






\bibitem{ILOG} Ikehata, M., Inverse obstacle scattering problems with a single incident wave and the logarithmic differential
               of the indicator function in the enclosure method, Inverse Problems, {\bf 27}(2011) 085006(23pp).





\bibitem{IE3} Ikehata, M.,
              The enclosure method for inverse obstacle scattering problems with dynamical data
              over a finite time interval: II.  Obstacles with a dissipative boundary or
              finite refractive index and back-scattering data, Inverse Problems, {\bf 28}(2012) 045010(29pp).



\bibitem{IE4} Ikehata, M.,
              An inverse acoustic scattering problem inside a cavity with dynamical back-scattering data,
              Inverse Problems, {\bf 28}(2012) 095016(24pp).





\bibitem{IE5} Ikehata, M.,
              The enclosure method for inverse obstacle scattering problems with dynamical data over a finite time interval: III.
              Sound-soft obstacle and bistatic data,
              Inverse Problems, {\bf 29}(2013) 085013(35pp).





\bibitem{IITOU} Ikehata, M. and Itou, H., Extracting the support function of a cavity in an isotropic elastic
                body from a single set of boundary data, Inverse Problems, {\bf 25}(2009) 105005(21pp).




\bibitem{IK2} Ikehata, M. and Kawashita, M.,
              On the reconstruction of inclusions in a heat conductive body
              from dynamical boundary data over a finite time interval,
              Inverse Problems, {\bf 26}(2010) 095004(15pp).







\bibitem{INS} Ikehata, M., Niemi, E. and Siltanen, S.,
              Inverse obstacle scattering with limited-aperture data,
              Inverse Problems and Imaging, {\bf 6}(2012), No.1, 77-94.




\bibitem{IO}  Ikehata, M. and Ohe, T., 
              Numerical method for finding the convex hull of polygonal cavities using enclosure method, 
              Inverse Problems, {\bf 18}(2002), 111-124. 




\bibitem{LP}  Lax, P. D. and Phillips, R. S.,  The scattering of sound waves by an obstacle,
              Comm. Pure Appl. Math., {\bf 30}(1977), 195-233.





\bibitem{M}   Majda, A., High frequency asymptotics for the scattering matrix and the inverse problem of
              acoustic scattering, Comm. Pure and Appl. Math., {\bf 29}(1976), 261-291.






\end{thebibliography}
\end{document}